\documentclass[12pt]{article} 
\usepackage{amssymb}          
\usepackage{latexsym}         

\font\csc=cmcsc10 scaled\magstep1

\font\teneufm=eufm10
\font\seveneufm=eufm7
\font\fiveeufm=eufm5
\newfam\eufmfam
\textfont\eufmfam=\teneufm
\scriptfont\eufmfam=\seveneufm
\scriptscriptfont\eufmfam=\fiveeufm

\newcommand{\slth}{\widehat{\mbox{\twelveeufm sl}}_2} 
\newcommand{\slN}{\mbox{\twelveeufm sl}_N} 
\font\twelveeufm=eufm10 scaled\magstep1
\font\seventeeneufm=eufm10 scaled\magstep3   
   
\newcommand{\slNBig}{\mbox{\seventeeneufm sl}_N} 
\newcommand{\Z}{{\mathbb Z}} 
\newcommand{\C}{{\mathbb C}} 

\newcommand{\F}{{\cal F}}
\newcommand{\At}{\tilde{A}}
\newcommand{\cL}{{\cal L}}
\newcommand{\ve}{\varepsilon}
\newcommand{\la}{\lambda}
\newcommand{\Wb}{\overline{W}}
\newcommand{\dz}{\underline{dz}}
\newcommand{\dw}{\underline{dw}}
\newcommand{\dbr}[1]{{\langle\!\langle #1 \rangle\!\rangle}} 
\newcommand{\Hc}{{\cal H}}
\newcommand{\G}[3]{{{G}_{x^{#1}}(#2,#3)}}
\newcommand{\ma}[1]{[\![ #1 ]\!]}
\newcommand{\dsp}{\displaystyle}
\newcommand{\qi}[1]{{ [  #1 ]_x}}
\newcommand{\dva}{${\cal V}_{x,r}(A^{(2)}_2)$ }
\newcommand{\dvaa}{${\cal V}_{x,r}(A^{(2)}_2)$}

\newcommand{\nn}{\nonumber}
\newcommand{\eqref}[1]{(\ref{#1})}
\newcommand{\bea}{\begin{eqnarray}}
\newcommand{\ena}{\end{eqnarray}}
\newcommand{\be}{\begin{eqnarray*}}
\newcommand{\en}{\end{eqnarray*}}
\newcommand{\lb}[1]{\label{#1}}

\newcommand{\res}{{\mathop{\rm res}}}
\newcommand{\id}{{\rm id}}
\newcommand{\tr}{{\rm tr}}
\newcommand{\bra}[1]{\langle #1 |}        
\newcommand{\ket}[1]{{| #1 \rangle}}      
\newcommand{\rs}{{r^*}}
\newcommand{\BW}[5]{\Bigl({#1\atop#3}\ {#2\atop#4}\Bigl|#5\Bigr)}
\newcommand{\Remark}{\medskip \noindent {\it Remark.}\quad}
\newcommand{\example}{\medskip \noindent {\it Example.}\quad}
\newcommand{\qed}{\hfill \fbox{}\medskip}
\newcommand{\proof}{\medskip\noindent{\it Proof.}\quad }

\newtheorem{thm}{Theorem}[section]
\newtheorem{conj}[thm]{Conjecture}
\newtheorem{prop}[thm]{Proposition}
\newtheorem{lem}[thm]{Lemma}
\newtheorem{dfn}[thm]{Definition}
\newcommand{\ignore}[1]{}

\begin{document}
\font\csc=cmcsc10 scaled\magstep1

{\baselineskip=14pt
 \rightline{
 \vbox{\today\hfill
       \hbox{HU-IAS/K-7}\\
       \phantom{a}\hfill\hbox{DPSU-99-2}
}}}

\vskip 11mm
\begin{center}
{\Large\bf Free Field Approach to the Dilute $A_L$ Models } 
\end{center}
\vskip11mm
\begin{center}
{\csc Yuji Hara}$\,^{1}$,
{\csc Michio Jimbo}$\,^{2}$,
{\csc Hitoshi Konno}$\,^{3}$,\\
{\csc Satoru Odake}$\,^{4}$ and
{\csc Jun'ichi Shiraishi}$\,^{5}$
\\ 
{\baselineskip=15pt
\it\vskip.35in 
\setcounter{footnote}{0}\renewcommand{\thefootnote}{\arabic{footnote}}
\footnote{e-mail address : snowy@gokutan.c.u-tokyo.ac.jp}
Institute of Physics, Graduate School of Arts and Sciences, \\
University of Tokyo, Tokyo 153-8902, Japan\\
\vskip.1in 
\footnote{e-mail address : jimbo@kusm.kyoto-u.ac.jp}
Division of Mathematics, Graduate School of Science,\\
Kyoto University, Kyoto 606-8502, Japan\\
\vskip.1in 
\footnote{e-mail address : konno@mis.hiroshima-u.ac.jp}
Department of Mathematics, Faculty of Integrated Arts and Sciences,\\
Hiroshima University, Higashi-Hiroshima 739-8521,  Japan\\
\vskip.1in 
\footnote{e-mail address : odake@azusa.shinshu-u.ac.jp}
Department of Physics, Faculty of Science \\
Shinshu University, Matsumoto 390-8621, Japan\\
\vskip.1in 
\footnote{e-mail address : shiraish@momo.issp.u-tokyo.ac.jp}
Institute for Solid State Physics, \\
University of Tokyo, Tokyo 106-0032, Japan \\
}
\end{center}

\vskip5mm

\begin{abstract}
We construct a free field realization of vertex operators 
of the dilute $A_L$ models along with the Felder complex. 
For $L=3$, we also study an $E_8$ structure 
in terms of the deformed Virasoro currents.
\end{abstract}

\vskip20mm
math.QA/9902150

\newpage

\setcounter{footnote}{0}
\renewcommand{\thefootnote}{\arabic{footnote})}

\setcounter{section}{0}
\setcounter{equation}{0}
\section{Introduction} \lb{sec:1}

The dilute $A_L$ model \cite{WNS92,WPSN94} is an integrable 
lattice model obtained by an RSOS restriction of 
the face model of type $A^{(2)}_2$ \cite{Ku91}. 
It possesses several intriguing features. 
Among others, we are interested in the following points.  
\begin{itemize}
\item[i)] At criticality, the dilute $A_L$ model in regime $2^+$ 
is described by conformal field theory (CFT) which belongs to 
the Virasoro minimal unitary series with the central charge 
${c=1-\frac{6}{L(L+1)}}$. 
The Andrews-Baxter-Forrester (ABF) model in regime III 
is known to be a different off-critical lattice model 
having the same critical behavior. 
While the ABF model corresponds to the ($1,3$)-perturbation 
of the minimal unitary CFT, 
the dilute $A_L$ model corresponds to the 
($1,2$)-perturbation of the same CFT. 

\item[ii)] In the particular case of $L=3$, 
the model falls within the same universality class 
as the two-dimensional Ising model in a magnetic field. 
The elliptic nome in the Boltzmann weights plays the role of a magnetic 
field, as opposed to the usual role as a temperature-like variable. 
In the field theory limit, the scattering process of particles 
exhibits an $E_8$ structure. 
\end{itemize}

In \cite{LP96}, bosonization of the ABF model in regime III was achieved. 
By `bosonization' we mean a free field realization 
of vertex operators (VO's) of the model as formulated in \cite{Fo94}.   
It was also found that the deformed Virasoro algebra (DVA) 
proposed earlier in \cite{qVir} arises naturally, 
in such a way that the VO's play the role of deformed 
chiral primary fields for DVA.  
(More specifically, 
the VO's of type I and type II correspond to the simplest 
primary fields $\phi_{21}$ and $\phi_{12}$, respectively. 
Analogs of general $\phi_{mn}$ are obtained by a
fusion procedure \cite{Kadei}. 
For an interpretation of the VO's as intertwiners 
for elliptic algebras, see \cite{JKOS2}.)
The BRST resolution of Fock spaces which singles out 
irreducible representations of the Virasoro algebra 
carries over to the deformed version as well. 
Thus the work \cite{LP96} presents 
an off-critical lattice version of the ($1,3$)-perturbation of 
the minimal unitary CFT in the free field picture. 
In this paper we study a similar problem for the ($1,2$)-perturbed CFT,  
by bosonizing the dilute $A_L$ model.  
For this purpose we adapt the construction of \cite{Fo94,LP96,MW96} 
to the case of $A^{(2)}_2$.  

In the trigonometric limit, 
bosonization of VO's has been given in \cite{Matsuno,JgMi97}  
using representation theory of the quantum affine algebra $U_q(A^{(2)}_2)$. 
In principle, we are to follow its elliptic analog 
on the basis of the face type elliptic algebra \cite{Konno,JKOS2}. 
Because of some technical difficulties in dealing with the latter, 
we take here a more pedestrian way 
and solve the exchange relations directly  
to obtain the bosonic realizations of VO's. 
In the course we use the elliptic Drinfeld currents 
obtained by a `dressing' procedure \cite{JKOS2}.  

In view of the points i), ii) above, it is natural to expect that 
a different deformation of the Virasoro algebra 
arises from the dilute $A_L$ models. 
Such a deformation was found in \cite{BL97} 
through bosonization of the $A^{(2)}_2$ affine Toda field theory.   
(See \cite{FR97,BP98} for more general deformed $W$ algebras
including the case $A^{(2)}_{2l}$.)
To make distinction from the original DVA of \cite{qVir},   
we use the symbol \dva to denote the DVA of \cite{BL97}. 
In the original case, the generating function of DVA 
(hereafter referred to as the `DVA current') 
can also be obtained from `fusion' of the VO's 
\cite{Luk96,JKM96,JS97}. 
In the same way, we reproduce the current of \dva 
by taking residues of products of the bosonized VO's
in the present case.

For the dilute $A_L$ models in regime $2^+$, 
the space of states of the corner transfer matrix 
is an analog of the minimal unitary representation \cite{WPSN94} for \dvaa. 
In order to obtain them from the Fock spaces, 
we consider a Felder type resolution 
using the elliptic currents of type $A^{(2)}_2$ as screening currents. 
Unlike the case of the ABF models, 
the BRST charges are not simply a power of screening operators. 
Such a complication seems to be common in 
the higher rank situations \cite{FJMOP}. 
We prove the nilpotency of BRST charges with the help of  
the Feigin-Odesskii algebra \cite{FJMOP}. 
Assuming a cohomological property of the resulting complex, 
we write down integral expressions for the two-point 
local height probabilities (LHPs) 
and traces of general product of the VO's. 

In the case of the dilute $A_3$ model,   
it is of some interest to see how Zamolodchikov's 
$E_8$ structure of scattering process \cite{Zam89} 
looks like in the free field picture. 
What plays the role of the Zamolodchikov-Faddeev operator 
for creation/annihilation of bound states is the DVA current \cite{Luk96}.  
Specializing to $L=3$, we introduce eight kinds of DVA currents by suitably 
fusing the elementary current of \dvaa. 
We then find a curious similarity between 
the operator product expansions of these currents and 
the so-called $T$-system of $E_8$ type 
for the transfer matrix \cite{KNS94,Su98}. 
At present we do not understand its proper meaning. 

This paper is organized as follows.
In section 2, we review the definition of the dilute $A_L$ model 
and give a brief description of the vertex operator approach. 
In section 3, we give bosonization of the VO's. 
We derive the current for \dva from them and 
state a conjecture for the Kac determinant formula for \dvaa. 
We also present a Felder type BRST complex of the Fock spaces.  
Section 4 is devoted to an application of the 
bosonization to the calculation of the LHP. 
In section 5, regarding the \dva 
current as the ZF operator for the particles in the $A_3$ model, 
we discuss 
a similarity between the $T$-system 
and the `bootstrap' of the ZF operators.  
Appendix A is a summary of the operator product 
expansion formulae for the elliptic currents and VO's. 
In appendix B, we prove the nilpotency of the BRST charges. 
Appendix C is devoted to an exposition of 
the fusion properties of the deformed $W$ currents for $A^{(1)}_{N-1}$, 
which should be compared with those for the DVA currents 
associated with the dilute $A_3$ model in section 5. 
 
\setcounter{section}{1}
\setcounter{equation}{0}
\section{The dilute $A_L$ models} \lb{sec:2}

\subsection{Boltzmann weights}\lb{subsec:2.1}

Throughout this paper we fix a positive integer $L\ge 3$. 
In  the dilute $A_L$ model, the local fluctuation variables 
$a,b,\cdots$ take one of the $L$ states $1,2,\cdots, L$, and 
those on neighboring lattice sites are subject to the condition $a-b=0,\pm 1$. 
The Boltzmann weights can be found in \cite{WPSN94}, eq.(3.1).
For our purpose it is convenient to use the 
parametrization given in Appendix A of \cite{WPSN94},  
which is suitable in the `low-temperature' regime. 
With some change of notation we recall the formula below.

Let $x=e^{-2\pi \lambda/\varepsilon}$, $r=\pi/(2\lambda)$ 
and $u=-u_{orig}/(2\lambda)$,  
where $\lambda$, $\varepsilon$ are the variables used in \cite{WPSN94} and 
$u_{orig}$ stands for `$u$' there. 
We shall restrict ourselves to the `regime $2^+$' defined by 
\bea
0<x<1, \quad 
r=2\,\frac{L+1}{L+2},
\quad
-\frac{3}{2}<u<0.
\lb{reg2}
\ena
Along with the variable $u$, we often use the multiplicative 
variable
\be
z=x^{2u}.
\en
Changing an overall scalar factor we put 
the Boltzmann weights in the form 
\be
&&W\BW{a}{b}{c}{d}{u}= \rho(u)\Wb\BW{a}{b}{c}{d}{u}, 
\en
where $\rho(u)$ will be specified in \eqref{rho} below.   
To give the formula for the $\Wb$ factors, let us set 
\bea
&&[u]=x^{\frac{u^2}{r}-u}\Theta_{x^{2r}}(x^{2u}),
\qquad 
[u]_+=x^{\frac{u^2}{r}-u}\Theta_{x^{2r}}(-x^{2u}),
\lb{[u]}
\ena
where 
\be
&&\Theta_p(z)=(z;p)_\infty(pz^{-1};p)_\infty(p;p)_\infty,
\\
&&(z;p_1,\cdots,p_k)_\infty=\prod_{n_1,\cdots,n_k=0}^\infty
(1-p_1^{n_1}\cdots p_k^{n_k}z).
\en
Then we have 
\begin{eqnarray*}
&&\Wb\BW{a\pm1}{a}{a}{a\mp1}{u}=1,
\\
&&\Wb\BW{a}{a\pm1}{a}{a\pm1}{u}=
\Wb\BW{a\pm1}{a\pm1}{a}{a}{u}=
-
\left(\frac{[\pm a+3/2]_+[\pm a-1/2]_+}{[\pm a+1/2]_+^2}\right)^{1/2}
\frac{[u]}{[1+u]},
\\
&&
\Wb\BW{a\pm1}{a}{a}{a}{u}=
\Wb\BW{a}{a}{a}{a\pm1}{u}=
\frac{[\pm a+1/2+u]_+}{[\pm a+1/2]_+}
\frac{[1]}{[1+u]},
\\
&&
\Wb\BW{a}{a\mp1}{a\pm1}{a}{u}=
\left(G_a^+G_a^-\right)^{1/2}\frac{[1/2+u]}{[3/2+u]}\frac{[u]}{[1+u]},
\\
&&
\Wb\BW{a}{a}{a\pm1}{a}{u}=
\Wb\BW{a}{a\pm1}{a}{a}{u}=
-\left(G_a^{\pm}\right)^{1/2}
\frac{[\pm a-1-u]_+}{[\pm a+1/2]_+}\frac{[1]}{[1+u]}\frac{[u]}{[3/2+u]},
\\
&&
\Wb\BW{a}{a\pm1}{a\pm1}{a}{u}=
\frac{[\pm 2a+1-u]}{[\pm 2a+1]}\frac{[1]}{[1+u]}
-
G_a^{\pm}\frac{[\pm 2a-1/2-u]}{[\pm 2a+1]}\frac{[u]}{[3/2+u]}\frac{[1]}{[1+u]},
\\
&&
\Wb\BW{a}{a}{a}{a}{u}=
\frac{[3+u]}{[3]}
\frac{[1]}{[1+u]}
\frac{[3/2-u]}{[3/2+u]}
+H_a\frac{[1]}{[3]}\frac{[u]}{[1+u]}.
\end{eqnarray*}
Here
\begin{eqnarray}
  G_a^{\pm}=\frac{S(a\pm1)}{S(a)},\quad
  S(a)=(-1)^{a}\frac{[2a]}{[a]_+},\quad
  H_a=G_a^+\frac{[a-5/2]_+}{[a+1/2]_+}+G_a^-\frac{[a+5/2]_+}{[a-1/2]_+}.
\end{eqnarray}
We choose $\rho(u)$ so that 
the partition function per site of the model equals to $1$. 
Explicitly it is given by \cite{WPSN94} 
\bea
z^{\frac{r-1}{r}}\rho(u)=\frac{\rho_+(u)}{\rho_+(-u)},\quad
\rho_+(u)=
\frac{(x^2z,x^3z, x^{2r+3}z, x^{2r+4}z;x^6,x^{2r})_\infty}
{(x^5z,x^6z,x^{2r}z,x^{2r+1}z;x^6,x^{2r})_\infty}, 
\lb{rho}
\ena
where $z=x^{2u}$, and we have introduced the notation 
\bea
(a_1,\cdots,a_n;p_1,\cdots,p_k)_\infty
=\prod_{j=1}^n(a_j;p_1,\cdots,p_k)_\infty.
\lb{mulpro}
\ena

Graphically we represent the Boltzmann weights as follows:
$$
\setlength{\unitlength}{1mm}
\begin{picture}(140,41)
\put(10,10){\makebox(50,20)[l]{${\displaystyle W\BW{a}{b}{c}{d}{u_1-u_2}\;=}$}}
\put(60,11){\vector(1,0){10}} \put(70,11){\line(1,0){8}}
\put(60,29){\vector(1,0){10}} \put(70,29){\line(1,0){8}}
\put(60,29){\vector(0,-1){10}} \put(60,19){\line(0,-1){8}}
\put(78,29){\vector(0,-1){10}} \put(78,19){\line(0,-1){8}}
\put(55,29){\makebox(5,5){$a$}}
\put(55,6){\makebox(5,5){$c$}}
\put(78,29){\makebox(5,5){$b$}}
\put(78,6){\makebox(5,5){$d$}}
\multiput(84.5,20)(-2,0){3}{\line(-1,0){1}}
\multiput(76.5,20)(-2,0){8}{\line(-1,0){1}}
\multiput(58.5,20)(-2,0){3}{\line(-1,0){1}}
\multiput(69,35.5)(0,-2){3}{\line(0,-1){1}}
\multiput(69,27.5)(0,-2){8}{\line(0,-1){1}}
\multiput(69,9.5)(0,-2){3}{\line(0,-1){1}}
\put(49.5,17.5){\makebox(5,5){$\lhd$}}
\put(66.5,1){\makebox(5,5){${\scriptstyle \bigtriangledown}$}}
\put(66.5,35.5){\makebox(5,5){$u_1$}}
\put(84.5,17.5){\makebox(5,5){$u_2$}}
\put(69,6){\makebox(5,5){$\ve_1$}}
\put(69,29){\makebox(5,5){$\ve'_1$}}
\put(55,15){\makebox(5,5){$\ve_2$}}
\put(78,15){\makebox(5,5){$\ve'_2$}}
\put(89.5,17.5){\makebox(5,5)[b]{,}}
\put(100,23.5){\makebox(20,5)[l]{$b=a+\ve'_1$,}}
\put(100,17.5){\makebox(20,5)[l]{$c=a+\ve_2$,}}
\put(100,11.5){\makebox(20,5)[l]{$d=b+\ve'_2=c+\ve_1$.}}
\end{picture}
$$

For definiteness we list below the basic properties of the Boltzmann weights. 
\begin{description}
\item[Yang-Baxter equation]
\be
&&\sum_{g}W\BW{a}{b}{g}{c}{u}W\BW{a}{g}{f}{e}{v}W\BW{g}{c}{e}{d}{u+v}
\\
&=\!\!&
\sum_{g}W\BW{a}{b}{f}{g}{u+v}W\BW{b}{c}{g}{d}{v}W\BW{f}{g}{e}{d}{u},
\en
\item[unitarity]
\[
\sum_{g}W\BW{a}{b}{g}{c}{u}
W\BW{a}{g}{d}{c}{-u}=\delta_{bd},
\]
\item[crossing symmetry]
\[
W\BW{b}{d}{a}{c}{-\frac{3}{2}-u}
=
\sqrt{\frac{S(a)S(d)}{S(b)S(c)}}
W\BW{a}{b}{c}{d}{u},
\]
\item[initial condition]
\[
W\BW{a}{b}{c}{d}{0}=\delta_{bc}.
\]
\end{description}

\subsection{Vertex operators}\lb{subsec:2.2}

Hereafter we assume that $L$ is odd. 
The model has ground states labeled by odd 
integers $l=1,3,\cdots,L-2$ \cite{WPSN94}.
They are characterized as configurations 
in which all heights take the same value $b$. 
If $L=4m\pm 1$, then the possible values are 
$b=l$ ($1\le l\le 2m-1$, $l$ : odd) or $b=l+1$ ($2m+1\le l\le L-2$, $l$ : odd).

Consider the corner transfer matrices $A(z),B(z),C(z),D(z)$ 
corresponding to the NW, SW, SE, NE quadrants, respectively. 
In the infinite volume limit, we have 
\[
C(z)=A(z)=z^{-\Hc},\quad B(z)=D(z)=\sqrt{S(k)}\,x^{3\Hc}z^{\Hc},  
\]
with $k$ denoting the value of the central height.
The operator $\Hc$ (the corner Hamiltonian) is independent of $z$.
We denote by $\cL_{l,k}$ the space of eigenstates of $\Hc$ 
in the sector where the central height is fixed to $k$ and  
the boundary heights are in the ground state $l$.
It was found in \cite{WPSN94} that  
the generating function of the spectrum of $\Hc$
coincides with the character of the Virasoro minimal unitary series.
Namely 
\bea
\tr_{\cL_{l,k}}(q^\Hc)=\chi_{l,k}(q),
\lb{char}
\ena
where 
\be
&&\chi_{l,k}(q)=
\frac{q^{\Delta_{l,k}-c/24}}{(q;q)_\infty}
\sum_{j\in \Z}\left(q^{L(L+1)j^2+((L+1)l-Lk)j}
-q^{L(L+1)j^2+((L+1)l+Lk)j+lk}\right),
\\
&&
c=1-\frac{6}{L(L+1)},\quad
\Delta_{l,k}=\frac{((L+1)l-Lk)^2-1}{4L(L+1)}.
\en

Consider next the half-infinite 
transfer matrix extending to infinity in the north. 
We denote it by 
\be
\Phi^{(k,k+\ve)}(z):\cL_{l,k+\ve}\rightarrow \cL_{l,k}
\qquad (\ve=0,\pm 1).
\en
Likewise we denote by 
\be
\Phi^{*(k+\ve,k)}(z^{-1}):\cL_{l,k}\rightarrow \cL_{l,k+\ve}
\qquad (\ve=0,\pm 1)
\en
the half-infinite transfer matrix extending to infinity in the west. 
We shall also write
\[
\Phi^{(k,k+\ve)}(z)=\Phi_\ve(z),
\quad
\Phi^{*(k+\ve,k)}(z)=\Phi^*_\ve(z), 
\]
and call them vertex operators (VO's) of type I.  
$$
\setlength{\unitlength}{1mm}
\begin{picture}(140,41)
\put(10,10){\makebox(20,20)[l]{$\Phi^{(a,b)}(z)\;=$}}
\put(35,33){\vector(0,-1){5}}
\put(35,28){\vector(0,-1){10}}
\put(35,18){\vector(0,-1){10}}
\put(35,8){\line(0,-1){4}}
\put(35,4){\vector(1,0){6}} \put(41,4){\line(1,0){4}}
\put(35,14){\vector(1,0){6}} \put(41,14){\line(1,0){4}}
\put(35,24){\vector(1,0){6}} \put(41,24){\line(1,0){4}}
\put(45,33){\vector(0,-1){5}}
\put(45,28){\vector(0,-1){10}}
\put(45,18){\vector(0,-1){10}}
\put(45,8){\line(0,-1){4}}
\put(35,4){\makebox(10,10){$u$}}
\put(35,14){\makebox(10,10){$u$}}
\put(35,24){\makebox(10,10){$u$}}
\multiput(40,33)(0,2){3}{\line(0,1){1}}
\put(30,-1){\makebox(5,5){$a$}}
\put(45,-1){\makebox(5,5){$b$}}
\put(50,17.5){\makebox(5,5)[b]{,}}
\put(60,10){\makebox(20,20)[l]{$\Phi^{*(a,b)}(z^{-1})\;=$}}
\put(96,25){\vector(1,0){5}}
\put(101,25){\vector(1,0){10}}
\put(111,25){\vector(1,0){10}}
\put(121,25){\line(1,0){4}}
\put(105,25){\vector(0,-1){6}} \put(105,19){\line(0,-1){4}}
\put(115,25){\vector(0,-1){6}} \put(115,19){\line(0,-1){4}}
\put(125,25){\vector(0,-1){6}} \put(125,19){\line(0,-1){4}}
\put(96,15){\vector(1,0){5}}
\put(101,15){\vector(1,0){10}}
\put(111,15){\vector(1,0){10}}
\put(121,15){\line(1,0){4}}
\put(95,15){\makebox(10,10){$u$}}
\put(105,15){\makebox(10,10){$u$}}
\put(115,15){\makebox(10,10){$u$}}
\multiput(96,20)(-2,0){3}{\line(-1,0){1}}
\put(125,10){\makebox(5,5){$a$}}
\put(125,25){\makebox(5,5){$b$}}
\put(130,17.5){\makebox(5,5)[b]{.}}
\end{picture}
$$

Intuitive graphical arguments 
based on the properties of the Boltzmann weights  
lead to the following formulas. 
For the details we refer the reader to \cite{JM,LP96}. 
\be
&&\Phi^{(a,c)}(z_2)\Phi^{(c,d)}(z_1)
=\sum_{g}
W\BW{a}{g}{c}{d}{u_1-u_2}
\Phi^{(a,g)}(z_1)\Phi^{(g,d)}(z_2)
\qquad (z_j=x^{2u_j}),\\
&&w^{\Hc}\Phi^{(a,b)}(z)w^{-\Hc}=\Phi^{(a,b)}(wz),
\qquad
\Phi^{*(b,a)}(z)=\sqrt{\frac{S(a)}{S(b)}}\Phi^{(b,a)}(x^{-3}z),
\\
&&\sum_{\ve}\Phi^{*(a,g)}(z)\Phi^{(g,b)}(z)=\delta_{ab},
\qquad
\Phi^{(a,b)}(z)\Phi^{*(b,c)}(z)=\delta_{ac}.
\en

As explained in \cite{Fo94}, 
multi-point local height probabilities 
are expressed as traces of VO's. 
Consider neighboring $n+1$ lattice sites in a row. 
Let $P_l(a_0,\cdots,a_n)$ denote the probability 
of finding these local variables to be $(a_0,\cdots,a_n)$. 
Then we have 
\bea
  \!\!\!\!\!&&P_l(a_0,\cdots,a_n)\nonumber\\
  \!\!\!\!\!&=\!\!&
  Z_{l}^{-1}S(a_n)\,\tr_{\cL_{l,a_n}}\Bigl(x^{6\Hc}
  \Phi^{*(a_n,a_{n-1})}(z)\cdots\Phi^{*(a_1,a_0)}(z)
  \Phi^{(a_0,a_1)}(z)\cdots\Phi^{(a_{n-1},a_n)}(z)\Bigr).
  \lb{corr}
\ena
Here the nomalization factor $Z_l$ is
\be
  Z_l=\sum_{k=1}^L S(k)\chi_{l,k}(x^6),
\en
which can be expressed in product of theta functions with conjugate 
modulus \cite{WPSN94}.
In the simplest case $n=0$, the one-point function 
$P_l(k)$ is given by 
$P_l(k)=Z_l^{-1}S(k)\chi_{l,k}(x^6)$.

\Remark 
We follow mostly the notation of \cite{MW96} but there are minor changes. 
The $\Phi(\zeta^{-1})$ in \cite{MW96} corresponds to 
$\Phi(z)$ in the present notation.
We have also reversed the orientation of edges of the Boltzmann weights.

\setcounter{section}{2}
\setcounter{equation}{0}
\section{Bosonization of vertex operators} \lb{sec:3}

\subsection{Bosons}\lb{subsec:3.1}

In this section we present a bosonic realization of vertex operators.
The working closely follows \cite{LP96,MW96}.
We set
\[
  [n]_x=\frac{x^n-x^{-n}}{x-x^{-1}},
\]
and introduce the oscillators $\alpha_n$ ($n\neq 0$) and 
$P,Q$ satisfying the commutation relations
\begin{eqnarray}
  &&
  [\alpha_n,\alpha_m]=
  \frac{[n]_x\left([2n]_x-[n]_x\right)}{n}
  \frac{[rn]_x}{[(r-1)n]_x}\delta_{n+m,0},\\
  &&[P,iQ]=1.\nonumber
\end{eqnarray}
Notice that $[2n]_x-[n]_x=[3n]_x[n/2]_x/[3n/2]_x$.
We shall also use
\be
  &&\alpha'_n=(-1)^n\frac{[(r-1)n]_x}{[rn]_x}\alpha_n.
\en
We denote by
\be
  &&\F_{l,k}=\C[\alpha_{-1},\alpha_{-2},\cdots]\ket{l,k}
\en
the Fock space generated by 
$$
  \ket{l,k}=e^{p_{l,k}iQ}\ket{0,0},\quad
  P\ket{l,k}=p_{l,k}\ket{l,k},
$$
where $p_{l,k}$ is 
\begin{equation}
  p_{l,k}=-\frac{l}{2}\sqrt{\frac{r}{r-1}}+k\sqrt{\frac{r-1}{r}}
  =-l\sqrt{\frac{L+1}{2L}}+k\sqrt{\frac{L}{2(L+1)}}.
  \label{plk}
\end{equation}
(Recall that $r=2(L+1)/(L+2)$.) 
These Fock spaces are graded by 
\be
&&d=\sum_{n=1}^\infty\frac{n^2}{[n]_x([2n]_x-[n]_x)}\frac{[(r-1)n]_x}{[rn]_x}
\alpha_{-n}\alpha_n
+\frac{1}{2}P^2-\frac{1}{24},
\en
which satisfies $[d,\alpha_n]=-n\alpha_n$, $[d,iQ]=P$ and 
$d\ket{l,k}=(\Delta_{l,k}-c/24)\ket{l,k}$. 
For later use, we define 
operators $\hat{l},\hat{k}:\F_{l,k}\rightarrow\F_{l,k}$ by 
\[
  \hat{l}\,|_{\F_{l,k}}=l\times\id_{\F_{l,k}},\quad
  \hat{k}\,|_{\F_{l,k}}=k\times\id_{\F_{l,k}}.
\]

\subsection{Vertex operators}\lb{subsec:3.2}

In the works \cite{Matsuno,JgMi97},  
a bosonic realization of the level-one representation of 
the quantum affine algebra $U_q(A^{(2)}_2)$   
and associated vertex operators have been obtained.
We shall consider their elliptic counterparts.
 
The elliptic version of the Drinfeld currents 
are constructed from the trigonometric ones by a
`dressing' procedure described in \cite{JKOS2}. 
Applying it to the present case of $U_q(A^{(2)}_2)$,  
we obtain 
\bea
  &&x_+(z):\F_{l,k}\rightarrow \F_{l-2,k},\quad
  x_-(z):\F_{l,k}\rightarrow \F_{l,k-1},\nn\\
  &&x_+(z)=\,
  :\exp\Biggl(-\sum_{n\neq 0}\frac{\alpha_n}{[n]_x}z^{-n}\Biggr):\times 
  e^{\sqrt{\frac{r}{r-1}}iQ}z^{\sqrt{\frac{r}{r-1}}P+\frac{r}{2(r-1)}},
  \lb{x+}\\
  &&x_-(z)=\,
  :\exp\Biggl(\sum_{n\neq 0}\frac{\alpha'_n}{[n]_x}z^{-n}\Biggr):\times
  e^{-\sqrt{\frac{r-1}{r}}iQ}z^{-\sqrt{\frac{r-1}{r}}P+\frac{r-1}{2r}}.
  \lb{x-}
\ena

The elliptic version of VO's (of type I and type II) are defined  
in terms of their trigonometric ones and a `twistor' given by 
an infinite product of the universal $R$ matrix \cite{JKOS1}.  
They satisfy the commutation relations of the type 
\eqref{com1}-\eqref{com3} in the next subsection. 
As we do not know how to evaluate the twistor in the bosonic realization,   
we have solved the relations \eqref{com1}-\eqref{com3} directly 
for $\Phi_\ve(z),\Psi^*_\ve(z)$. 
We obtain the following. 
\medskip

\noindent {\bf Type I}:
\bea
&&\Phi_\ve(z):\F_{l,k}\rightarrow \F_{l,k-\ve},\nonumber\\
  \Phi_-(z)&\!\!=\!\!&
  :\exp\Biggl(-\sum_{n\neq 0}\frac{\alpha'_n}{[2n]_x-[n]_x}z^{-n}\Biggr):\times
  e^{\sqrt{\frac{r-1}{r}}iQ}z^{\sqrt{\frac{r-1}{r}}P+\frac{r-1}{2r}},
  \lb{Phi-}
  \\
  \Phi_0(z)&\!\!=\!\!&
  x^{\frac{1-r}{2r}}\oint_{C_0}\dz_1\Phi_-(z)x_-(z_1)
  \frac{1}{\sqrt{[\hat{k}+1/2]_+[\hat{k}-1/2]_+}}
  \frac{[u-u_1+\hat{k}]_+}{[u-u_1+1/2]},
  \lb{Phi0}\\
  \Phi_+(z)&\!\!=\!\!&
  x^{\frac{1-r}{r}}\oint\!\!\oint_{C_+}\dz_1\dz_2\Phi_-(z)x_-(z_1)x_-(z_2)\nn\\
  &&\times\sqrt{\frac{S(\hat{k}-1)}{S(\hat{k})}}
  \frac{1}{[\hat{k}-1/2]_+[2\hat{k}-2]}
  \frac{[u-u_1+2\hat{k}-3/2]}{[u-u_1+1/2]}
  \frac{[u_1-u_2+\hat{k}]_+}{[u_1-u_2+1/2]}.
  \lb{Phi+}
\ena

\noindent {\bf Type II}:
\bea
&&\Psi^*_\ve(z):\F_{l,k}\rightarrow \F_{l-2\ve,k},\nonumber\\
  \!\!\!\!\!\!\!
  \Psi^*_-(z)&\!\!=\!\!&
  :\exp\Biggl(\sum_{n\neq 0}\frac{\alpha_n}{[2n]_x-[n]_x}z^{-n}\Biggr):\times
  e^{-\sqrt{\frac{r}{r-1}}iQ}z^{-\sqrt{\frac{r}{r-1}}P+\frac{r}{2(r-1)}},
  \lb{Psi-}\\
  \!\!\!\!\!\!\!
  \Psi^*_0(z)&\!\!=\!\!&
  ix^{\frac{r}{2(r-1)}}\oint_{C^*_0}\dz_1\Psi^*_-(z)x_+(z_1)
  \frac{1}{\sqrt{[(\hat{l}+1)/2]^*_+[(\hat{l}-1)/2]^*_+}}
  \frac{[u-u_1-\hat{l}/2]^*_+}{[u-u_1-1/2]^*},
  \lb{Psi0}\\
  \!\!\!\!\!\!\!
  \Psi^*_+(z)&\!\!=\!\!&
  x^{\frac{r}{r-1}}\oint\!\!\oint_{C^*_+}\dz_1\dz_2
  \Psi^*_-(z)x_+(z_1)x_+(z_2)\nn\\
  \!\!\!\!\!\!\!
  &&\!\!\!\!\!\times\sqrt{\frac{S^*(\hat{l}/2-1)}{S^*(\hat{l}/2)}}
  \frac{1}{[(\hat{l}-1)/2]^*_+[\hat{l}-2]^*}
  \frac{[u-u_1-\hat{l}+3/2]^*}{[u-u_1-1/2]^*}
  \frac{[u_1-u_2-\hat{l}/2]^*_+}{[u_1-u_2-1/2]^*}.
  \lb{Psi+}
\ena
Here $z=x^{2u}$, $z_j=x^{2u_j}$, $\dz_j=dz_j/(2\pi iz_j)$ and
\be
  [u]^*=x^{\frac{u^2}{r-1}-u}\Theta_{x^{2r-2}}(x^{2u}),\quad 
  [u]^*_+=x^{\frac{u^2}{r-1}-u}\Theta_{x^{2r-2}}(-x^{2u}),\quad
  S^*(a)=(-1)^a\frac{[2a]^*}{[a]^*_+}.
\en
The poles of the integrand of \eqref{Phi0}--\eqref{Psi+}
and the integration contours are listed in the following table
($n=0,1,2,\cdots$). 
For example, $C_0$ is a simple closed contour 
that encircles $x^{1+2rn}z$ ($n\ge 0$) but not $x^{-1-2rn}z$ ($n\ge 0$).  
$$
\begin{tabular}{|c|c|c|}
\hline
&inside&outside\\
\hline
$C_0$&$z_1=x^{1+2rn}z$&$z_1=x^{-1-2rn}z$\\
\hline
$C_+$&$z_1=x^{1+2rn}z$&$z_1=x^{-1-2rn}z$\\
&$z_2=x^{1+2rn}z_1$&$z_2=x^{-1-2rn}z,x^{-1-2rn}z_1,x^{2-2r(n+1)}z_1$\\
\hline
$C^*_0$&$z_1=x^{-1+2(r-1)n}z$&$z_1=x^{1-2(r-1)n}z$\\
\hline
$C^*_+$&$z_1=x^{-1+2(r-1)n}z$&$z_1=x^{1-2(r-1)n}z$\\
&$z_2=x^{-1+2(r-1)n}z_1$&
$z_2=x^{1-2(r-1)n}z,x^{1-2(r-1)n}z_1,x^{-2-2(r-1)(n+1)}z_1$\\
\hline
\end{tabular}
$$

\subsection{Commutation relations and inversion identities}\lb{subsec:3.3}

The VO's given above satisfy the following commutation relations.
\begin{eqnarray}
  \!\!\!\!\!\!\!\!\!\!\!\!\!\!\!
  &&\Phi_{\ve_2}(z_2)\Phi_{\ve_1}(z_1)
  =\!\!\!\!\!\sum_{\ve_1',\ve_2'\atop \ve_1'+\ve_2'=\ve_1+\ve_2}\!\!\!\!\!
  W\BW{\hat{k}}{\hat{k}+\ve_1'}{\hat{k}+\ve_2}{\hat{k}+\ve_1+\ve_2}{u_1-u_2}
  \Phi_{\ve'_1}(z_1)\Phi_{\ve'_2}(z_2),
  \lb{com1}\\
  \!\!\!\!\!\!\!\!\!\!\!\!\!\!\!
  &&\Psi^*_{\ve_1}(z_1)\Psi^*_{\ve_2}(z_2)
  =\!\!\!\!\!\sum_{\ve_1',\ve_2'\atop \ve_1'+\ve_2'=\ve_1+\ve_2}\!\!\!\!\!
  W^*\BW{\hat{l}/2}{\hat{l}/2+\ve_1}{\hat{l}/2+\ve_2'}
  {\hat{l}/2+\ve_1+\ve_2}{u_1-u_2}
  \Psi^*_{\ve'_2}(z_2)\Psi^*_{\ve'_1}(z_1),
  \lb{com2}\\
  \!\!\!\!\!\!\!\!\!\!\!\!\!\!\!
  &&\Phi_{\ve_2}(z_2)\Psi^*_{\ve_1}(z_1)=
  \tau(u_1-u_2)\Psi^*_{\ve_1}(z_1)\Phi_{\ve_2}(z_2).
\lb{com3}
\end{eqnarray}
Here we have set (for $z=x^{2u}$)  
\bea
  &&W^*\BW{a}{b}{c}{d}{u}=
  \Wb\BW{a}{b}{c}{d}{u}\Biggl|_{r\rightarrow r-1}
  \times \rho^*(u),
\lb{W*}\\
  &&z^{-\frac{r}{r-1}}\rho^*(u)=
  \frac{\rho^*_+(u)}{\rho^*_+(-u)},\quad
  \rho^*_+(u)=
  \frac{(x^3z,x^4z,x^{2r}z,x^{2r+1}z;x^6,x^{2r-2})_\infty}
  {(z,xz,x^{2r+3}z,x^{2r+4}z;x^6,x^{2r-2})_\infty},
  \lb{rho*}\\
  &&
  \tau(u)=z\frac{\Theta_{x^6}(-xz^{-1})\Theta_{x^6}(-x^2z^{-1})}
  {\Theta_{x^6}(-xz)\Theta_{x^6}(-x^2z)}.
\lb{tau}
\ena
Note that
\[
  \rho^*(u)=-\rho(u)\Bigl|_{r\rightarrow r-1}\times 
  z\frac{\Theta_{x^6}(xz^{-1})\Theta_{x^6}(x^2z^{-1})}
        {\Theta_{x^6}(xz)\Theta_{x^6}(x^2z)}.
\]
We do not present the tedious but straightforward 
verification of \eqref{com1}-\eqref{com3}.

For the description of correlation functions we need also the `dual' VO's. 
Define
\begin{eqnarray}
  &&\Phi^*_\ve(z)=g\sqrt{S(\hat{k})}^{-1}\Phi_{-\ve}(x^{-3}z)
  \sqrt{S(\hat{k})},\\
  &&\Psi_\ve(z)=g^{*-1}\sqrt{S^*(\hat{l}/2)}\Psi^*_{-\ve}(x^{-3}z)
  \sqrt{S^*(\hat{l}/2)}^{-1},
\end{eqnarray}
where
\begin{eqnarray*}
  &&g^{-1}=\frac{(x;x^{2r})_\infty}
  {(x^2;x^{2r})^2_\infty(x^{2r-1};x^{2r})_\infty(x^{2r};x^{2r})^4_\infty}
  \frac{(x^5,x^6,x^{2r},x^{2r+1};x^6,x^{2r})_\infty}
  {(x^2,x^3,x^{2r+3},x^{2r+4};x^6,x^{2r})_\infty},\\
  &&g^*=\frac{(x^{-1};x^{2r-2})_\infty}
  {(x^{-2};x^{2r-2})^2_\infty(x^{2r-1};x^{2r-2})_\infty
  (x^{2r-2};x^{2r-2})^5_\infty}
  \frac{(x^3,x^4,x^{2r},x^{2r+1};x^6,x^{2r-2})_\infty}
  {(x,x^6,x^{2r+3},x^{2r+4};x^6,x^{2r-2})_\infty}.
\end{eqnarray*}
Then we have 
\ignore{
\bea
  \!\!\!\!\!\!\!\!\!\!\!\!\!\!\!\!\!\!\!\!
  &&\Phi_{\ve_2}(z)\Phi^*_{\ve_1}(z)=\delta_{\ve_1,\ve_2}\times\id,\quad
  \sum_{\ve}\Phi^*_{\ve}(z)\Phi_{\ve}(z)=\id,
  \lb{inv1}\\
  \!\!\!\!\!\!\!\!\!\!\!\!\!\!\!\!\!\!\!\!
  &&\Psi_{\ve_1}(z_1)\Psi^*_{\ve_2}(z_2)
  =\frac{\delta_{\ve_1,\ve_2}}{1-z_1/z_2}+\cdots,\;
  \sum_{\ve}\Psi^*_{\ve}(z_2)\Psi_{\ve}(z_1)=\frac{1}{1-z_1/z_2}+\cdots,
  \;(z_1\rightarrow z_2).
\lb{inv2}
\ena
}
\bea
  \!\!\!\!\!\!\!\!\!\!
  &&\Phi_{\ve_2}(z)\Phi^*_{\ve_1}(z)=\delta_{\ve_1,\ve_2}\times\id,\quad
  \Psi_{\ve_1}(z_1)\Psi^*_{\ve_2}(z_2)
  =\frac{\delta_{\ve_1,\ve_2}}{1-z_1/z_2}+\cdots,\quad
  (z_1\rightarrow z_2),
  \lb{inv1}\\
  \!\!\!\!\!\!\!\!\!\!
  &&\sum_{\ve}\Phi^*_{\ve}(z)\Phi_{\ve}(z)=\id,\qquad\quad
  \sum_{\ve}\Psi^*_{\ve}(z_2)\Psi_{\ve}(z_1)=\frac{1}{1-z_1/z_2}+\cdots,
  \quad (z_1\rightarrow z_2).
\lb{inv2}
\ena

\subsection{Deformed Virasoro algebra}\lb{subsec:3.4}

Brazhnikov and Lukyanov \cite{BL97} pointed out that 
one can associate to the algebra $A^{(2)}_2$ 
a deformed Virasoro algebra (DVA) which is different from the 
one found in \cite{qVir}. 
The original DVA of \cite{qVir}, associated with $A^{(1)}_1$,  
arises also as a `fusion' of VO's \cite{Luk96}.  
Let us discuss this point in the present case of $A^{(2)}_2$.

Let
\begin{eqnarray}
  &&\Lambda_{\pm}(z)=\,:\exp\Biggl(\pm\sum_{n\neq0}
  \lambda_n(x^{\pm3/2}z)^{-n}\Biggr):\times x^{\pm 2\sqrt{r(r-1)}P},
\nonumber\\
  &&\Lambda_{0}(z)=-\frac{[r-1/2]_x}{[1/2]_x}
  :\exp\Biggl(\sum_{n\neq0}\lambda_n(x^{-n/2}-x^{n/2})z^{-n}\Biggr):,
  \lb{Tz}\\
%
  &&T(z)=\Lambda_+(z)+\Lambda_0(z)+\Lambda_-(z),\nonumber
\end{eqnarray}
where
\begin{eqnarray*}
  &&\lambda_n=(-1)^n(x-x^{-1})\frac{[(r-1)n]_x}{[2n]_x-[n]_x}\alpha_n
  =(x-x^{-1})\frac{[rn]_x}{[2n]_x-[n]_x}\alpha'_n,\\
  &&[\lambda_n,\lambda_m]=(x-x^{-1})^2
  \frac{1}{n}\frac{[n]_x[rn]_x[(r-1)n]_x}{[2n]_x-[n]_x}\delta_{m+n,0}.
\end{eqnarray*}
Then $T(z)$ is obtained from VO's by fusing them,
\begin{eqnarray}
  \!\!\!\!\!\!\!\!\!\!
  &&\Phi_{\ve_2}(x^{r+3/2}z')\Phi^*_{\ve_1}(x^{-r+3/2}z)\nn\\
  \!\!\!\!\!\!\!\!\!\!&=\!\!&
  \left(1-\frac{z}{z'}\right)(-1)^{\ve_1+1}
  \delta_{\ve_1,\ve_2}T(z)\cdot x^{1-r}
  \frac{(x,x^6,x^{5-2r},x^{6-2r};x^6)_\infty}
       {(x^3,x^4,x^{2-2r},x^{3-2r};x^6)_\infty}
  +\cdots,\quad(z'\rightarrow z).
  \lb{DVA}
\end{eqnarray}

The $T(z)$ satisfies the DVA of \cite{BL97}
\begin{eqnarray}
  &&f\left(\frac{z_2}{z_1}\right)T(z_1)T(z_2)
  -f\left(\frac{z_1}{z_2}\right)T(z_2)T(z_1)\nonumber\\
  &=\!\!&(x-x^{-1})\frac{[r+1/2]_x[r]_x[r-1]_x[r-3/2]_x}
  {[1/2]_x[3/2]_x}\left(\delta\Bigl(x^3\frac{z_2}{z_1}\Bigr)
  -\delta\Bigl(x^{-3}\frac{z_2}{z_1}\Bigr)\right)
  \label{BLDVA}\\
  &&+(x-x^{-1})\frac{[r]_x[r-1/2]_x[r-1]_x}{[1/2]_x}
  \left(\delta\Bigl(x^2\frac{z_2}{z_1}\Bigr)T(xz_2)
  -\delta\Bigl(x^{-2}\frac{z_2}{z_1}\Bigr)T(x^{-1}z_2)\right),\nonumber
\end{eqnarray}
where $f(z)$ is
\begin{eqnarray}
  f(z)&\!\!=\!\!&
  \exp\Biggl(-\sum_{n>0}(x-x^{-1})^2
  \frac{1}{n}\frac{[n]_x[rn]_x[(r-1)n]_x}{[2n]_x-[n]_x}z^n\Biggr)\nonumber\\
  &\!\!=\!\!&
  \frac{1}{1-z}
  \frac{(x^{2-2r}z,x^{3-2r}z,x^4z,x^5z,x^{2r}z,x^{2r+1}z;x^6)_{\infty}}
       {(x^{5-2r}z,x^{6-2r}z,xz,x^2z,x^{2r+3}z,x^{2r+4}z;x^6)_{\infty}}.
  \label{fz}
\end{eqnarray}
The notation of \cite{BL97} is related to ours by 
$x_{BL}=x^{3/2}$, $b/Q=r$, $1/(Qb)=1-r$, $g(z)=f(z)$, ${\bf V}(z)=T(z)$.
In what follows we call this algebra \dvaa.
The relation \eqref{BLDVA} is invariant under 
\begin{equation}
  r\mapsto 1-r,\quad x\mapsto x,\quad T(z)\mapsto -T(z). 
  \label{rto1-r}
\end{equation}
We remark that $\widetilde{T}(z)=-\Lambda_+(z)+\Lambda_0(z)-\Lambda_-(z)$ 
also satisfies (\ref{BLDVA}), which is obtained from type II VO's,
\be
  &&\Psi_{\ve_1}(x^{r-1+3/2}z')\Psi^*_{\ve_2}(x^{-(r-1)+3/2}z)\nn\\
  &=\!\!&
  \frac{1}{1-z'/z}(-1)^{\ve_1+1}\delta_{\ve_1,\ve_2}
  \widetilde{T}(-z)\cdot(-x^{-r})
  \frac{(x^2,x^3,x^{5-2r},x^{6-2r};x^6)_\infty}
       {(x^5,x^6,x^{2-2r},x^{3-2r};x^6)_\infty}
  +\cdots,\quad(z'\rightarrow z).
\en

Let us discuss some features of \dvaa. 
\medskip

\noindent {\bf Conformal limit}\quad 
In the conformal limit ($x=e^{\hbar}\rightarrow 1$, $r$ : fixed), 
\eqref{BLDVA} admits two limits \cite{BL97} related by (\ref{rto1-r}),
\begin{eqnarray}
  \!\!\!\!\!\!\!\!\!\!\!\!\!\!\!\!\!\!\!\!&&
  T(z)=3-2r+\hbar^2\Bigl(8r(r-1)z^2L(z)+\frac16r(r-1)(1-2r)+(2-r)^2\Bigr)
  +O(\hbar^4),
  \label{3-2r}\\
  \!\!\!\!\!\!\!\!\!\!\!\!\!\!\!\!\!\!\!\!&&
  T(z)=-1-2r+\hbar^2\Bigl(-8r(r-1)z^2\widetilde{L}(z)
  +\frac16r(r-1)(1-2r)-(1+r)^2\Bigr)+O(\hbar^4),
  \label{-1-2r}
\end{eqnarray}
where $L(z), \widetilde{L}(z)$ are the Virasoro currents with the central 
charges $c,\widetilde{c}$ respectively,
$$
  c=1-\frac{3(2-r)^2}{r(r-1)}=1-\frac{6}{L(L+1)},\quad
  \widetilde{c}=1-\frac{3(1+r)^2}{r(r-1)}.
$$

In the free boson realization \eqref{Tz}, $T(z)$ and $\widetilde{T}(z)$ 
have `natural' expansions \eqref{3-2r} and \eqref{-1-2r} respectively,
by the following identification:
\be
  &&\lambda_n=2\hbar\sqrt{r(r-1)}\sqrt{\frac{x^{rn}-x^{-rn}}{2\hbar rn}
  \frac{x^{(r-1)n}-x^{-(r-1)n}}{2\hbar(r-1)n}\frac{1}{x^n-1+x^{-n}}}a_n,\\
  &&P=a_0-\frac{2-r}{2\sqrt{r(r-1)}}
  =\widetilde{a}_0-\frac{r+1}{2\sqrt{r(r-1)}},\\
  &&L(z)=\,:\frac12\Bigl(\partial\phi(z)\Bigr)^2:+\frac{2-r}{2\sqrt{r(r-1)}}
  \partial^2\phi(z),\\
  &&\widetilde{L}(z)=\,:\frac12\Bigl(\partial\widetilde{\phi}(z)\Bigr)^2:
  +\frac{r+1}{2\sqrt{r(r-1)}}\partial^2\widetilde{\phi}(z),
\en
where $[a_n,a_m]=n\delta_{n+m,0}$,  
$\partial \phi(z)=\sum_{n\in\Z}a_nz^{-n-1}$ and
$\partial\widetilde{\phi}(z)=\partial\phi(z)
\Bigl|_{a_0\rightarrow\widetilde{a}_0}$.
On the other hand, $T(z)$ has an expansion of the form 
\eqref{-1-2r} with $P=\widetilde{a}_0-\frac{r+1}{2\sqrt{r(r-1)}}+
\frac{i\pi(2n+1)}{2\hbar\sqrt{r(r-1)}}$ ($n\in\Z$).
 
\medskip

\noindent {\bf Kac determinant}\quad 
Let $T(z)=\sum_{n\in\Z}T_nz^{-n}$, 
and let $U_\pm$ be the algebra generated by $\{T_n\}_{\pm n>0}$. 
As usual, the Verma module of highest weight $\lambda\in\C$ is defined 
as the free left $U_-$-module generated by 
a vector $\ket{\lambda}$ such that 
$T_n\ket{\lambda}=0$ ($n>0$) and $T_0\ket{\lambda}=\lambda\ket{\lambda}$.  
Likewise the right Verma module is defined by 
$\bra{\lambda}T_n=0$ ($n<0$), 
$\bra{\lambda}T_0=\lambda\bra{\lambda}$, and
$\langle\lambda|\lambda\rangle=1$.
At level $N$ there are $p(N)$ (the number of partition) independent 
states, $T_{-n_1}T_{-n_2}\cdots T_{-n_l}\ket{\lambda}$ 
($n_1\geq n_2\geq\cdots\geq n_l>0$, $\sum_{i=1}^ln_i=N$).
Let us number these states by the reverse lexicographic ordering for 
$(n_1,n_2,\cdots,n_l)$, i.e., 
$\ket{\lambda;N,1}=T_{-N}\ket{\lambda}$, 
$\ket{\lambda;N,2}=T_{-N+1}T_{-1}\ket{\lambda}$, $\cdots$, 
$\ket{\lambda;N,p(N)}=T_{-1}^N\ket{\lambda}$.
Similarly we define  
$\bra{\lambda;N,1}=\bra{\lambda}T_N$, 
$\bra{\lambda;N,2}=\bra{\lambda}T_1T_{N-1}$, $\cdots$,  
$\bra{\lambda;N,p(N)}=\bra{\lambda}T_1^N$.  

We conjecture that the Kac determinant at level $N$ is given by 
\bea
  &&\det\Bigl(\langle\lambda;N,i|\lambda;N,j\rangle\Bigr)_{1\leq i,j\leq p(N)}
  \nn\\
  &=\!\!&\prod_{l,k\geq 1 \atop lk\leq N}\left(
  (\lambda-\lambda_{l,k})(\lambda-\widetilde{\lambda}_{l,k})
  \frac{(x^{rl}-x^{-rl})(x^{(r-1)l}-x^{-(r-1)l})}{x^l-1+x^{-l}}
  \right)^{p(N-lk)},
\ena
where
\be
  \lambda_{l,k}&\!\!=\!\!&
  x^{-lr+2k(r-1)}+x^{lr-2k(r-1)}-\frac{[r-1/2]_x}{[1/2]_x},\\
  \widetilde{\lambda}_{l,k}&\!\!=\!\!&
  -x^{-2lr+k(r-1)}-x^{2lr-k(r-1)}-\frac{[r-1/2]_x}{[1/2]_x}.
\en
We remark that in the free boson realization 
$T_0\ket{l,k}=\lambda_{l,k}\ket{l,k}$ and 
$T_0e^{\alpha iQ}\ket{2l,k/2}=
\widetilde{\lambda}_{l,k}e^{\alpha iQ}\ket{2l,k/2}$,
where $\alpha=\frac{i\pi(2n+1)}{2\hbar\sqrt{r(r-1)}}$ ($n\in\Z$).

\subsection{Felder complex}\lb{subsec:3.5}

The Fock spaces $\F_{l,k}$ themselves do not give
a bosonic realization of the space of states 
${\cal L}_{l,k}$ of the corner Hamiltonian.  
For this we need a cohomological construction using 
an analog of the Felder complex \cite{Fel89}:
\begin{eqnarray}
&&
\cdots ~{\buildrel X_{-2} \over \longrightarrow }~
\F_{2L-l,k} ~{\buildrel X_{-1} \over \longrightarrow }~
\F_{l,k} ~{\buildrel X_0 \over \longrightarrow }~
\F_{-l,k} ~{\buildrel X_1 \over \longrightarrow }~
\F_{l-2L,k} ~{\buildrel X_2 \over \longrightarrow }~
\cdots, 
\lb{BRS}
\\
&&\quad X_jX_{j-1}=0.\nn
\end{eqnarray}
In the case of the algebra $A^{(1)}_1$, Lukyanov and Pugai constructed 
the coboundary map $X_j$ as a power of a single operator 
(\cite{LP96}, see also \cite{JLMP}). 
In our case, the formula for $X_j$ is a little more involved.   

Set 
\be
Q_1&\!\!=\!\!&\oint_{|z|=1}\dz\,x_+(z)\frac{[u+\hat{l}/2]^*}{[u+1/2]^*},
\\
Q_2^{(a)}&\!\!=\!\!&\oint\!\!\oint_{|z_1|=|z_2|=1}
\dz_1\dz_2\,x_+(z_1)x_+(z_2)
\frac{1}{[u_1+1/2]^*[u_2+1/2]^*}
\\
&&\times \frac{[u_1-u_2]^*}{[u_1-u_2+1]^*[u_1-u_2-1/2]^*}
f_2^{(a)}(u_1+\hat{l}/2,u_2+\hat{l}/2),
\en
where
\be
f^{(a)}_2(u_1,u_2)&\!\!=\!\!&
[2a+1]^*[a-1/2]^*[u_1-a]^*[u_2+a-1]^*[u_1-u_2+a-1/2]^*
\\
&&\!\!\!\!-[2a-1]^*[a+1/2]^*[u_1+a]^*[u_2-a-1]^*[u_1-u_2-a-1/2]^*.  
\en
These operators are mutually commutative (see Lemma \ref{lem3}). 
We define `BRST charges' $Q_l$ ($1\le l\le L-1$) as follows: 
\be
Q_l=\cases{Q_1Q_2^{(1)}\cdots Q_2^{(m)} & ($l=2m+1$),\cr
Q_2^{((L+1)/2-m)}\cdots Q_2^{((L-3)/2)}Q_2^{((L-1)/2)} & ($l=2m$).\cr}
\en
Note that $Q_L=Q_lQ_{L-l}$. 

We prove the following propositions in appendix \ref{app:2}
(Proposition \ref{prop:app3.1}, \ref{prop:app3.2}). 
\begin{prop}\lb{prop:3.51}
Suppose $l'$ is odd and $l'\equiv l~\bmod~L$ ($1\le l\le L-1$). 
On the space $\F_{l',k}$, $Q_l$ is expressed as
\bea
&&Q_l=\oint\!\cdots\!\oint_{|z_1|=\cdots=|z_l|=R}\dz_1\cdots\dz_l
\,x_+(z_1)\cdots x_+(z_l)\,H_l(u_1,\cdots,u_l),
\lb{Ql2}\\
&&
H_l(u_1,\cdots,u_l)=\pm\bar{h}_l(u_1,\cdots,u_l)\prod_{1\le i<j\le l}
\frac{[u_i-u_j]^*}{[u_i-u_j+1]^*[u_i-u_j-1/2]^*},
\lb{Hl}
\ena
where $\bar{h}_l(u_1,\cdots,u_l)$ is holomorphic, symmetric and satisfies
$\bar{h}_l(u_1+v,\cdots,u_l+v)=\bar{h}_l(u_1,\cdots,u_l)$. 
We have
\bea
&&H_l(u_1+r-1,\cdots,u_l)=H_l(u_1,\cdots,u_l),
\lb{Hl1}\\
&&H_l(u_1+\tau,\cdots,u_l)=H_l(u_1,\cdots,u_l)\,e^{-\pi i(l-1)/(r-1)},
\lb{Hl2}
\ena
where $\tau=\pi i/\log x$.
\end{prop}
Hence \eqref{Ql2} does not depend on $R>0$. 

\begin{prop}\lb{prop:3.52}
Under the same condition as above, we have 
\be
Q_lQ_{L-l}=0\qquad (1\le l\le L-1). 
\en
\end{prop}

Let us call $C_{l,k}$ the cochain complex \eqref{BRS}
defined by
\begin{eqnarray*}
X_{2j}=Q_l&:&\F_{l-2jL,k}\longrightarrow \F_{-l-2jL,k}, 
\\
X_{2j+1}=Q_{L-l}&:&\F_{-l-2jL,k}\longrightarrow \F_{l-2(j+1)L,k}. 
\end{eqnarray*}
In the conformal limit where $x\rightarrow 1$ and $z=x^{2u}$ kept fixed,  
this complex formally tends to Felder's complex \cite{Fel89} 
for the minimal unitary series. 
In view of this, it is natural to expect that 
\begin{equation}
H^j(C_{l,k})={\rm Ker} X_{j}/{\rm Im} X_{j-1}=0 \qquad (j\neq 0).  
\lb{coho}
\end{equation}
By Euler-Poincar\'e principle, the $0$-th cohomology $H^0(C_{l,k})$ 
has then the same character as the space of states 
${\cal L}_{l,k}$ (see \eqref{char}),
\[
\tr_{H^0(C_{l,k})}(q^d)=\tr_{{\cal L}_{l,k}}(q^{\Hc}).
\]
Henceforth we assume \eqref{coho} and make an identification 
\[
H^0(C_{l,k})={\cal L}_{l,k},\quad d=\Hc.
\]  

\begin{prop}\lb{prop:3.53}
Under the same assumption as in Proposition \ref{prop:3.51}, 
we have on $\F_{l',k}$	
\bea
&&[\Phi_{\ve}(z),Q_l]=0\quad (\ve=0,\pm),	
\lb{com4}\\
&&[T(z),Q_l]=0. \lb{com5}
\ena
\end{prop}
\medskip

\proof 
\eqref{com5} is a consequence of \eqref{com4} and \eqref{DVA}.
Let us show
\[
\Phi_-(z)Q_l=Q_l\Phi_-(z).
\]
The left (resp. right) hand side is well defined
if we choose $R\ll 1$ (resp. $R\gg 1$) in \eqref{Ql2}.
As meromorphic functions we have 
$\Phi_-(z)x_+(z_j)=x_+(z_j)\Phi_-(z)$,  
and the product has no poles. 
Since $Q_l$ does not depend on $R$, the conclusion follows. 

Next let us prove \eqref{com4} with $\ve=0$. 
The case $\ve=+$ can be shown similarly.
Dropping irrelevant constants we consider 
\[
\Phi_0'(z)=\oint_{C_0}\dz' \Phi_-(z)x_-(z')
\frac{[u-u'+\hat{k}]_+}{[u-u'+1/2]}.
\]
As meromorphic functions we have 
\[
x_-(z)x_+(z')=x_+(z')x_-(z)=\frac{z+z'}{(z+xz')(z+x^{-1}z')}
:x_-(z)x_+(z'):.
\]
We use the expression \eqref{Ql2} with $x^2|z|<R<x^{-2}|z|$.  
Taking into account the symmetry 
in the integration variables $z_1,\cdots,z_l$, we obtain 
\be
\!\!\!\!\!\!\!\!\!\!
&&[\Phi_0'(z),Q_l]=
l\oint\!\cdots
\!\oint_{|z_1|=\cdots=|z_l|=R}\dz_1\cdots\dz_l\,H_l(u_1,\cdots,u_l)
\\
\!\!\!\!\!\!\!\!\!\!
&&\qquad\qquad\qquad\times
\left(\res_{z'=-xz_1}+\res_{z'=-x^{-1}z_1}\right)
\Phi_-(z)x_-(z')x_+(z_1)\cdots x_+(z_l)
\frac{[u-u'+\hat{k}]_+}{[u-u'+1/2]}\dz'.
\en
By noting the identity 
\[
:x_+(z)x_-(-x^{-1}z):\,
=x^{-2r+1}:x_+(x^{2r-2}z)x_-(-x^{2r-1}z):, 
\]
we can rewrite the right hand side as follows:
\bea
\!\!\!\!\!&&l\oint\!\cdots\!\oint_{|z_2|=\cdots=|z_l|=R}\dz_2\cdots\dz_l
\Bigl(\oint_{C_1}\dz_1
A(z_1,z)x_+(z_2)\cdots x_+(z_l)H_l(u_1,\cdots,u_l)
\nn\\
\!\!\!\!\!&&\qquad\qquad\qquad\qquad\quad -
\oint_{C_2}\dz_1
A(x^{2r-2}z_1,z)x_+(z_2)\cdots x_+(z_l)H_l(u_1,\cdots,u_l)
\Bigr),
\lb{AA}
\ena
where
\be
A(z_1,z)&\!\!=\!\!&\res_{z'=-xz_1}\Phi_-(z)x_-(z')x_+(z_1)
\frac{[u-u'+\hat{k}]_+}{[u-u'+1/2]}\dz'\\
&\!\!=\!\!&\frac{\mbox{ holomorphic function }}
{(-x^2z_1/z,-x^{2r}z/z_1;x^{2r})_\infty}
:\Phi_-(z)x_-(-xz_1)x_+(z_1):.
\en
The contours for $z_1$ are ($n\ge 0$) 
$$
\begin{tabular}{|c|c|c|}
\hline
&inside&outside\\
\hline
$C_1$&$z_1=-x^{2r(n+1)}z$&$z_1=-x^{-2-2rn}z$\\
\hline
$C_2$&$z_1=-x^{2+2rn}z$&$z_1=-x^{-2r(n+1)}z$\\
\hline
\end{tabular}
$$
Moreover the product $:x_-(-xz_1)x_+(z_1):x_+(z_j)$ 
is holomorphic in $z_1$ for $|x^{2r}z_j|<|z_1|$. 
In view of the periodicity \eqref{Hl1}, 
the two terms of \eqref{AA} cancel out 
by shifting the contour $z_1\rightarrow x^{2r-2}z_1$. 
\qed

\setcounter{section}{3}
\setcounter{equation}{0}

\section{Local height probabilities} \lb{sec:4}
We present here a calculation of the local height probabilities (LHP) 
for the dilute $A_L$ models in the regime $2^+$.

\subsection{Two-point LHP}

We have already mentioned the result \eqref{char} 
about the one-point function.
As the next simplest case, let us 
consider the probability $ P_l(a-\ve,a)$ of finding two 
neighboring local height variables to be $a-\ve,\ a\ (\ve=0, \pm)$.
\bea
P_l(a-\ve,a)&\!\!=\!\!&\frac{1}{Z_l}\ S(a)\ 
\tr_{\cL_{l,a}}\Bigl(x^{6\Hc}\Phi^{*}_\ve(z)\Phi_\ve(z)\Bigr)\nonumber\\
&\!\!=\!\!&\frac{1}{Z_l}\ g \sqrt{S(a)S(a-\ve)}\ \tr_{\cL_{l,a}}
\Bigl(x^{6\Hc}\Phi_{-\ve}(x^{-3}z)\Phi_\ve(z)\Bigr)
\label{tpf}.
\ena
Note that $P_l(a-\ve,a)$ is independent of $z$. 
{}From \eqref{tpf} and the property of the type I VO \eqref{inv2}, 
we have the following relations.
\be
\sum_{\ve=\pm1,0} \ P_l(a-\ve,a)=\frac{S(a)\ \chi_{l,a}(x^6)}{Z_l},\quad
P_l(a-\ve,a)=P_l(a,a-\ve),\quad
P_l(0,1)=0.
\en

The evaluation of the trace yields the following expressions.
\be
&&P_l(a-1,a)=-\frac{S(a-1)\ x^{\frac{1-r}{r}}}{[a-\frac{1}{2}]_+[2a-2]}
\oint\oint_{C_+(1)} {\dw}_1\ \dw_2\ {\cal I}(w_1,w_2)\\
&&\qquad\qquad\qquad\qquad\qquad\qquad\qquad\qquad\times 
\frac{[v_1-\frac{1}{2}][v_1-2a+\frac{3}{2}]
[v_1-v_2+a]_+}
{[v_1+\frac{1}{2}][v_1-\frac{1}{2}][v_1-v_2+\frac{1}{2}]},\\
&&P_l(a,a)=\frac{S(a)\ x^{\frac{1-r}{r}} }{[a+\frac{1}{2}]_+[a-\frac{1}{2}]_+}
\oint_{C_0(x^{-3})}\dw_1\oint_{C_0(1)} \dw_2\ {\cal I}(w_1,w_2)
\frac{[v_1-a+\frac{3}{2}]_+[v_2-a]_+}
{[v_1+1][v_2-\frac{1}{2}]},\\
&&P_l(a+1,a)=-\frac{S(a)\ x^{\frac{1-r}{r}} }{[a+\frac{1}{2}]_+[2a]}
\oint\oint_{C_+(x^{-3})} \dw_1\ \dw_2\ {\cal I}(w_1,w_2)\\
&&\qquad\qquad\qquad\qquad\qquad\qquad\qquad\qquad\times 
\frac{[v_1-2a+1][v_1-v_2+a+1]_+[v_2+\frac{1}{2}]}
{[v_1+1][v_1-v_2+\frac{1}{2}][v_2-\frac{1}{2}]}.
\en
Here $w_i=x^{2v_i}$ ($i=1,2$) and 
\be
  \!\!\!\!\!
  &&{\cal I}(w_1,w_2)\\
  &=\!\!&\tr_{\cL_{l,a}}\Bigl(x^{6\Hc}\Phi^{}_-(x^{-3})x_-(w_1)
  \Phi_-(1)x_-(w_2)\Bigr)\\
  \!\!\!\!\!
  &=\!\!&{\cal O}_{l,k}(w_1,w_2) \frac{(x^5,x^6;x^6,x^{2r})_\infty^2}
  {(x^{2r+3},x^{2r+4};x^6,x^{2r})_\infty^2}
  \frac{(x^{2r+2}{w_1},x^{2r+2}/w_1,x^{2r-1}w_2,x^{2r+2}/w_2;
        x^3,x^{2r})_\infty}
  {(x^{4}{w_1},x^{4}/w_1,x^{}w_2,x^{4}/w_2;x^3,x^{2r})_\infty}\\
  \!\!\!\!\!
  &&\times(x^6,w_2/w_1,x^6w_1/w_2;x^6)_\infty
  \frac{\G{6}{x^{2r-1}}{w_2/w_1}\G{6}{x^{2}}{w_2/w_1}}{
  \G{6}{x^{2r-2}}{w_2/w_1}\G{6}{x}{w_2/w_1}},
\en
where ${\cal O}_{l,k}(w_1,w_2)$ is the zero-mode contribution
\be
  {\cal O}_{l,k}(w_1,w_2)&\!\!=\!\!&
  \left(x^{6}\right)^{\Delta_{l,a}-\frac{c}{24}}
  (x^3w_1w_2)^{\frac{l(L+1)-aL}{2(L+1)}+\frac{L}{4(L+1)}}
  \sum_{j\in \Z}(x^3w_1w_2)^{-Lj}\\
  &&\times
  \left(\left(x^{6}\right)^{ L(L+1)j^2-(l(L+1)-aL)j}-
  \left(x^{6}\right)^{ L(L+1)j^2+(l(L+1)+aL)j +la}
  (x^3 w_1w_2)^{-l}\right),
\en 
and 
\be
\G{6}{A}{z}=(x^6A;x^6,x^{2r})_\infty^2(Az;x^6,x^{2r})_\infty 
(x^6A/z;x^6,x^{2r})_\infty.
\en
The contours $C_+(1)$, $C_0(x^{-3})\cup C_0(1)$, $C_+(x^{-3})$ are chosen
as follows $(n,m\geq 0)$;
For all the contours,
the poles $w_1=x^{4+3m+2rn}$, $w_2=x^{4+3m+2rn}$, $x^{4+6m+2r(n+1)}w_1$,
$x^{1+6(m+1)+2rn}w_1$ are inside and the poles $w_1=x^{-4-3m-2rn}$, 
$w_2=x^{-1-3m-2rn}$, $x^{2-6m-2r(n+1)}w_1$, $x^{-1-6m-2rn}w_1$ are outside.
In addition,  
$$ 
\begin{tabular}{|c|c|c|}
\hline
&inside&outside\\
\hline
$C_+(1)$&$w_1=x^{-1+2r(n+1)}$&$w_1=x^{-1-2rn}$\\
&$w_2=x^{1+2rn}w_1$&$w_2=x^{-1-2rn},x^{-1-2rn}w_1,x^{2-2r(n+1)}w_1$\\
\hline
$C_0(x^{-3})\cup C_0(1)$&$w_1=x^{-2+2rn}$&$w_1=x^{-4-2rn}$\\
&$w_2=x^{1+2rn}$&$w_2=x^{-1-2rn}$\\
\hline
$C_+(x^{-3})$&$w_1=x^{-2+2rn}$&$w_1=x^{-4-2rn}$\\
&$w_2=x^{1+2rn},x^{1+2rn}w_1$&$w_2=x^{1-2r(n+1)},x^{-1-2rn}w_1,
x^{2-2r(n+1)}w_1$\\
\hline
\end{tabular}
$$

\subsection{General case}

Integral representation of the $N$-point 
correlation functions can be derived in a similar manner.  
It is written in terms of the traces of the type I vertex operators as
in \eqref{corr}:
\bea
&&
Z_{l}^{-1}S(k)\;\tr_{\cL_{l,k}}
\Bigl(\Phi^*_{\ve_1}(x^6z_1)\cdots\Phi^*_{\ve_N}(x^6z_N)
\Phi_{\ve_N}(x^6z_N)\cdots\Phi_{\ve_1}(x^6z_1)x^{6\Hc}\Bigr).
\lb{N-pt}
\ena
Here we give only the integral formula for the traces 
over the Fock module in a general situation 
\bea
\tr_{\F_{l,k}}\Bigl(\Phi_{\ve_1}(z_1)\cdots\Phi_{\ve_N}(z_N)x^{6\Hc}\Bigr).
\lb{trace}
\ena
We assume $\sum_{t=1}^N\ve_t=0$. Otherwise \eqref{trace} vanishes.  

First we prepare several functions.  
\be
&&F(z)=\frac{(x^{5+2r}z;x^6,x^{2r})}{(x^7z;x^6,x^{2r})},\qquad
G(z)=\frac{F(z)}{F(xz)F(x^{-1}z)},\\
&&H(z)=\frac{(x^{8}z,x^{9}z,x^{9+2r}z,x^{10+2r}z;x^6,x^6,x^{2r})}
{(x^{11}z,x^{12}z,x^{6+2r}z,x^{7+2r}z;x^6,x^6,x^{2r})}.
\en
Define $h_{\varepsilon_n}(z_n,\{w_{n,i}\},\hat{k})$ 
($\varepsilon_n=0,+$) by normal-ordering the integrand of 
\eqref{Phi0},\eqref{Phi+},  
\be
&&\Phi_{\varepsilon_n}(z_n)=\oint\Bigr(\prod_{i\in I(\ve)}\dw_{n,i}\Bigl)
\;:\Phi_-(z_n)\prod_{i\in I(\ve)}x_-(w_{n,i}):\;
h_{\varepsilon_n}(z_n,\{w_{n,i}\},\hat{k}),\\
&&
I(0)=\{1\},\;I(+)=\{1,2\}.
\en
Explicitly we have 
\be
&&
h_0(z,w,k)=\frac{x^{(2u-2v+k+1/2)(k-1/2)/r^2-k+1/2}
(xz^2)^{\frac{1-r}{2r}}}{\sqrt{[k+1/2]_+[k-1/2]_+}}
\frac{\displaystyle{(-x^{2k}z/w,-x^{2r-2k}w/z;x^{2r})_\infty}}
{\displaystyle{(xw/z,xz/w;x^{2r})_\infty}},\\[0.2cm]
&&
h_+(z,w_1,w_2,k)=\\
&&\qquad\sqrt{\frac{S(k-1)}{S(k)}}
\frac{x^{\{4(k-1)(u-v_1+k-1/2)+(k-1/2)(2v_1-2v_2+k+1/2)\}/r^2-3k+5/2}
(xz^2/w_1)^{\frac{1-r}{r}}}{[k-1/2]_+[2k-2]}\\[0.2cm]
&&
\qquad\times
\frac{\displaystyle{(x^{2r-1}w_2/z,x^{4k-3}z/w_1,
x^{2r-4k+3}w_1/z;x^{2r})_\infty}}
{\displaystyle{(xw_2/z,xw_1/z,xz/w_1;x^{2r})_\infty}}\\[0.2cm]
&&
\qquad\times
\frac{\displaystyle{(x^{2}w_2/w_1,-x^{2k}w_1/w_2,
-x^{2r-2k}w_2/w_1;x^{2r})_\infty}}
{\displaystyle{(xw_2/w_1,x^{2r-2}w_2/w_1,xw_1/w_2;x^{2r})_\infty}}
\Bigl(1-\frac{w_1}{w_2}\Bigr).
\en
We use the symbol $\dbr{A(z)B(w)}$ to denote the 
normal ordering factors 
\[
A(z)B(w)= \dbr{A(z)B(w)}:A(z)B(w):.
\]
(See the list in Appendix \ref{app:1}.)

With this notation we have 
\bea
  &&\tr_{\F_{l,k}}
  \Bigl(\Phi_{\ve_1}(z_1)\cdots\Phi_{\ve_N}(z_N)x^{6\Hc}\Bigr)\nn\\
  &=\!\!&
  \oint\cdots\oint\prod_{{1\leq m\leq N}\atop{\ve_m\ne -}}
  \Bigl(\prod_{j\in I(\ve_m)}\dw_{m,j}\Bigr)\;h_{\ve_m}
  (z_m,\{w_{m,j}\},k+\sum_{t=1}^{m}\ve_t)\nn\\
  &&\times
  \prod_{1\leq m<n\leq N}\dbr{\Phi_-(z_m)\Phi_-(z_n)}
  \prod_{{1\leq m<n\leq N}\atop{i\in I(\ve_n),j\in I(\ve_m)}}
  \dbr{x_-(w_{m,j})x_-(w_{n,i})}\\
  &&\times
  \prod_{{1\leq n<m\leq N}\atop{j\in I(\ve_m)}}
  \dbr{\Phi_-(z_n)x_-(w_{m,j})}
  \prod_{{1\leq m<n\leq N}\atop{j\in I(\ve_m)}}
  \dbr{x_-(w_{m,j})\Phi_-(z_n)}\nn \\
  &&\times
  \tr_{\F_{l,k}}\Bigl(:\Phi_-(z_1)\cdots\Phi_-(z_N)
  \prod_{{1\leq m\leq N}\atop{j\in I(\ve_m)}}x_-(w_{m,j}):x^{6\Hc}\Bigr),\nn
\ena
where 
\bea
&&
\tr_{\F_{l,k}}\Bigl(:\Phi_-(z_1)\cdots\Phi_-(z_N)
\prod_{{1\leq m\leq N}\atop{j\in I(\ve_m)}}x_-(w_{m,j}):x^{6\Hc}\Bigr)\nn\\
&=\!\!&
\prod_{1\leq m,n\leq N}H(z_n/z_m)
\prod_{{1\leq m,n\leq N}\atop{i\in I(\ve_n),j\in I(\ve_m)}}
G(w_{n,i}/w_{m,j})\prod_{{1\leq n,m\leq N}\atop{j\in I(\ve_m)}}
F(z_n/w_{m,j})F(w_{m,j}/z_n)
\lb{N-tr}\\
&&
\times\;\frac{x^{6(p_{l,k}^2/4-1/24)}}{(x^6;x^6)_\infty}\;
\left(\Bigl(\prod_{{1\leq m\leq N}
\atop{j\in I(\ve_m)}}w_{m,j}\Bigr)
\Big/\Bigl(\prod_{1\leq n\leq N} z_n\Bigr)
\right)^{\sqrt{ \frac{\scriptstyle{r-2}}{\scriptstyle{r}} }p_{l,k}} 
\Bigl(\prod_{{1\leq m\leq N}\atop{j\in I(\ve_m)}}w_{m,j}
\prod_{1\leq n\leq N} z_n
\Bigr)^{ \frac{\scriptstyle{r-1}}{\scriptstyle{2r}} }\nn
\ena
where $p_{l,k}$ is given in \eqref{plk}.
The following are 
the list of poles of the integrand as functions of $w_{m,j}$. 
The contour for $\dw_{m,j}$ encircles only those denoted \lq inside\rq 
($a,b\in\Z_{\ge 0}$) :
$$
\begin{tabular}{|c|c|c|}
\hline
&inside&outside\\
\hline
$h_0(z_m,w_{m,1},k)$&$w_{m,1}=x^{1+2rb}z_m$&$x^{-1-2rb}z_m$\\
\hline
$h_+(z_m,w_{m,1},w_{m,2},k)$&$w_{m,1}=x^{1+2rb}z_m$&
$w_{m,1}=x^{-1-2rb}z_m$\\
&$w_{m,2}=x^{1+2rb}w_{m,1}$&$w_{m,2}=x^{-1-2rb}w_{m,1},x^{-1-2rb}z_{m}$\\
\hline
$\dbr{x_-(w_{m,j})\Phi_-(z_n)}$&$w_{m,j}=x^{1+2rb}z_n$&\\
\hline
$\dbr{\Phi_-(z_n)x_-(w_{m,j})}$&&$w_{m,j}=x^{-1-2rb}z_n$\\
\hline
$\dbr{x_-(w_{m,j})x_-(w_{n,i})}$&
$w_{m,j}=x^{1+2rb}w_{n,i}$&\\
&$w_{m,j}=x^{-2+2r(1+b)}w_{n,i}$&\\
\hline
$G(w_{n,i}/w_{m,j})$&$w_{m,j}=x^{7+2rb+6a}w_{n,i}$&\\
&$w_{m,j}=x^{2r(1+b)+6(1+a)}w_{n,i}$&\\
&$w_{m,j}=x^{4+2r(1+b)+6a}w_{n,i}$&\\
\hline
$F(z_n/w_{m,j})$&$w_{m,j}=x^{7+2rb+6a}z_n$&$w_{m,j}=x^{-7-2rb-6a}z_n$\\
\hline
\end{tabular}
$$
The formula for the $N$-point correlation function \eqref{N-pt} 
can be obtained through specializing \eqref{N-tr} and noting 
\[
H(x^3z)H(z)=\frac{1}{F(x^2z)F(xz)}.
\]
Since the result is lengthy we do not present it here. 

\setcounter{section}{4}
\setcounter{equation}{0}
\section{Discussion}

As was discussed in the main text, the 
DVA for the 
dilute $A_L$ model (which we have denoted by \dvaa) 
exactly coincides with the
one found by Brazhnikov and Lukyanov \cite{BL97}.
In the paper \cite{BL97}, \dva with $|x|=1$ was treated as
the Zamolodchikov-Faddeev (ZF) algebra 
for the Bullough-Dodd model ($A^{(2)}_2$ Toda field theory).
We regard \dva with $0<x<1,r=2(L+1)/(L+2)$ 
as the ZF algebra for the 
dilute $A_L$ model (restricted face model), and
apply the idea of bootstrap method to study 
the fusion of the \dva current $T(z)$. 

The two-dimensional Ising model at the critical temperature $T=T_c$ is 
described by the $c=1/2$ minimal CFT. 
Perturbing it by a magnetic field 
while keeping the same temperature ($T=T_c$),  
an off-critical integrable model is obtained \cite{Zam89}.
A fascinating feature of this theory is that 
the Lie algebra $E_8$ appears as a hidden symmetry; 
one can check that the integrals of motion $P_s$ appear 
at the exponents of $E_8$, $s=1,7,11,13,17,19,\cdots$, 
the bootstrap program closes within eight particles,
the mass ratios are given by the Perron-Frobenius vector for
the incidence matrix of $E_8$, and so on. 
Further discussions of the model 
as the $\phi_{1,2}$-perturbation of the $c=1/2$ 
CFT can be found in \cite{Smi91}. 
It is argued that the dilute $A_3$ model is 
in the universality class of the magnetic-perturbed 
Ising model \cite{WNS92}. 
As in the case of the ABF model \cite{LP96}, 
our free field realization for the dilute $A_3$ model  
properly reduces to that of the $c=1/2$ CFT,  
including the VO's, \dva and the Felder complex.  
Our description of the the dilute $A_3$ model, therefore, 
provides a lattice analogue of the 
$\phi_{1,2}$-perturbation of the $c=1/2$ CFT.

In this section, we study
an $E_8$-structure arising from \dva 
for the dilute $A_3$ model ($r=8/5$).
We construct eight fused DVA currents $T^{(a)}(u)$ ($a=1,2,\cdots,8$) 
from the fundamental \dva current $T(z)$ using a bootstrap procedure. 
We show that these fused currents obey a set of relations which resembles 
the so called level-two restricted 
$T$-system of type $E^{(1)}_8$ \cite{KNS94}. 

The $T$-system of type $E^{(1)}_8$ \cite{KNS94} is written as 
\begin{eqnarray}
T^{(a)}_m(u-{1\over 20})T^{(a)}_m(u+{1\over 20})=
T^{(a)}_{m-1}(u)T^{(a)}_{m+1}(u)+
g^{(a)}_m(u) \prod_{b\sim a}T^{(b)}_m(u), \label{t-s}
\end{eqnarray}
where $T^{(a)}_m(u)$ $(a=1,2,\cdots,8)$ denotes 
the eigenvalues of the transfer matrix, the symbol $b\sim a$ 
means that $b$ and $a$ are adjacent nodes in the Dynkin diagram of $E_8$, 
and $g^{(a)}(u)$'s are some functions. 
(We have rescaled the $u$ variable of \cite{KNS94} 
to fit the present notation.)
If the face model is restricted, then 
we have the truncation $1\leq m \leq \ell$. 
The integer $\ell$ is called level. 
If $\ell=2$, $T^{(a)}_{2}(u)$ becomes proportional to the identity, 
and \eqref{t-s} reduces to 
\begin{eqnarray}
T^{(a)}_1(u-{1\over 20})T^{(a)}_1(u+{1\over 20})=
\phi^{(a)}(u)+
g^{(a)}_1(u) \prod_{b\sim a}T^{(b)}_1(u), \label{t-s1}
\end{eqnarray}
with some functions $\phi^{(a)}(u)$. 
This is called the level-two restricted $T$-system of type $E^{(1)}_8$.
In the paper \cite{Su98}, 
the $T$-system in \eqref{t-s1} is realized in terms of the  
`quantum' transfer matrix for the dilute $A_3$ model.
\bigskip

Now we come back to the deformed Virasoro algebra \dva
for the dilute $A_3$ model. 
Before going into the technical details, 
let us roughly state the type of formulas 
we find for \dva with $r=8/5$.   
\begin{dfn}
Define the eight DVA currents  $T^{(a)}(u)$ ($a=1,2,\cdots,8$) 
corresponding to the simple roots of $E_8$ by
\begin{eqnarray*}
&&T^{(a)}(u)= T_{\overline{a}}(u) \qquad\qquad a=1,2,3,4,5,\\
&&T^{(6)}(u)= f_{1\overline{3}}(u_2-u_1)T_{1}(u_1)T_{\overline{3}}(u_2)
\Biggr|_{u_1=u+{11\over 20} \atop
         u_2=u-{3\over 20}},\\
&&T^{(7)}(u)= T_{2}(u),\\
&&T^{(8)}(u)=
f_{1\overline{2}}(u_2-u_1)T_{1}(u_1)T_{\overline{2}}(u_2)
\Biggr|_{u_1=u+{9\over 20} \atop
         u_2=u-{4\over 20}}.
\end{eqnarray*}
Here the fused DVA currents $T_n(u), T_{\overline{n}}(u)$
and the structure function $f_{m\overline{n}}(u)$
are defined in Definition \ref{fused dva} below.
\end{dfn}

\begin{prop}\label{thm}
The following relations hold:
\begin{eqnarray}
\!\!\!\!\!\!\!\!\!\!\!\!\!\!\!&& f^{(1)}(u_1,u_2) T^{(1)}(u_1)T^{(1)}(u_2)
\Biggr|_{u_1=u+{1\over 20} \atop u_2=u-{1\over 20}}
= T^{(2)}(u),\label{ts-1}\\
\!\!\!\!\!\!\!\!\!\!\!\!\!\!\!&& f^{(2)}(u_1,u_2) T^{(2)}(u_1)T^{(2)}(u_2)
\Biggr|_{u_1=u+{1\over 20} \atop u_2=u-{1\over 20}}
= g^{(2)}(u_1,u_2) T^{(1)}(u_1)T^{(3)}(u_2)
\Biggr|_{u_1=u \atop u_2=u },\label{ts-2}\\
\!\!\!\!\!\!\!\!\!\!\!\!\!\!\!&& f^{(3)}(u_1,u_2) T^{(3)}(u_1)T^{(3)}(u_2)
\Biggr|_{u_1=u+{1\over 20} \atop u_2=u-{1\over 20}}
= g^{(3)}(u_1,u_2) T^{(2)}(u_1)T^{(4)}(u_2)
\Biggr|_{u_1=u \atop u_2=u },\label{ts-3}\\
\!\!\!\!\!\!\!\!\!\!\!\!\!\!\!&& f^{(4)}(u_1,u_2) T^{(4)}(u_1)T^{(4)}(u_2)
\Biggr|_{u_1=u+{1\over 20} \atop u_2=u-{1\over 20}}
= g^{(4)}(u_1,u_2) T^{(3)}(u_1)T^{(5)}(u_2)
\Biggr|_{u_1=u \atop u_2=u },\label{ts-4}\\
\!\!\!\!\!\!\!\!\!\!\!\!\!\!\!&& f^{(5)}(u_1,u_2) T^{(5)}(u_1)T^{(5)}(u_2)
\Biggr|_{u_1=u+{1\over 20} \atop u_2=u-{1\over 20}}
= g^{(5)}(u_1,u_2,u_3) T^{(4)}(u_1)T^{(6)}(u_2)T^{(8)}(u_3)
\Biggr|_{u_1=u \atop{ u_2=u \atop u_3=u} },\label{ts-5}\\
\!\!\!\!\!\!\!\!\!\!\!\!\!\!\!&& f^{(6)}(u_1,u_2) T^{(6)}(u_1)T^{(6)}(u_2)
\Biggr|_{u_1=u+{1\over 20} \atop u_2=u-{1\over 20}}
= g^{(6)}(u_1,u_2) T^{(5)}(u_1)T^{(7)}(u_2)
\Biggr|_{u_1=u \atop u_2=u },\label{ts-6}\\
\!\!\!\!\!\!\!\!\!\!\!\!\!\!\!&& f^{(7)}(u_1,u_2) T^{(7)}(u_1)T^{(7)}(u_2)
\Biggr|_{u_1=u+{1\over 20} \atop u_2=u-{1\over 20}}
= T^{(6)}(u),\label{ts-7}\\
\!\!\!\!\!\!\!\!\!\!\!\!\!\!\!&& f^{(8)}(u_1,u_2) T^{(8)}(u_1)T^{(8)}(u_2)
\Biggr|_{u_1=u+{1\over 20} \atop u_2=u-{1\over 20}}
= T^{(5)}(u),\label{ts-8}
\end{eqnarray}
with appropriate functions $f^{(a)},g^{(a)}$ (see Definition \ref{fg}).
Both sides are regarded as operators on the cohomology $H^0(C_{l,k})$.  
Likewise we have the relations
\begin{eqnarray}
&&\ma{u_2-u_1-{3\over 2}} f^{(a)}(u_1,u_2) 
T^{(a)}(u_1)T^{(a)}(u_2)
\Biggr|_{u_1=u+{1\over 20}-{r\over 2} 
\atop u_2=u-{1\over 20}+{r\over 2}}
= c^{(a)} {\rm id},\label{ts-id}
\end{eqnarray}
for $a=1,2,\cdots,8$, where $c^{(a)}$'s are some constants and 
the symbol $\ma{~}$ is defined below in \eqref{ma}.
\end{prop}
We notice that \eqref{ts-1}--\eqref{ts-id} for the fused \dva currents 
look very similar to the $T$-system \eqref{t-s1} arising 
from the analytic Bethe ansatz. 
There is, however, an obvious discrepancy between them; 
while the $T$-system \eqref{t-s1} comprises two terms, 
\eqref{ts-1}--\eqref{ts-id} consists of only one term. 
More precisely, the right hand side of \eqref{ts-1}--\eqref{ts-8} 
corresponds to the second term in \eqref{t-s1},
whereas that of \eqref{ts-id} corresponds to the first term. 
In the left hand side,  
the spectral parameters $u_1-u_2$ of \eqref{ts-1}--\eqref{ts-8} and 
\eqref{ts-id} differ by $r$. 
Such a shift by $r$ is irrelevant in the $T$-system \eqref{t-s1}, 
because the transfer matrix eigenvalues $T^{(a)}_m(u)$ 
(with appropriate normalization) are periodic, $T^{(a)}_m(u+r)=T^{(a)}_m(u)$.  
This is a reflection of the quasi-periodicity of the Boltzmann weights. 
On the other hand, the \dva currents $T^{(a)}(u)$ 
are {\it by no means} doubly quasi-periodic;
we have $T^{(a)}(u+\pi i/\log x)=T^{(a)}(u)$ but 
$T^{(a)}(u+r)\neq T^{(a)}(u)$. 

We have not understood yet the reason why we have 
such similarities between the $T$-system for the Bethe ansatz and the 
exchange relations for the DVA. 
For the purpose of comparison, 
we summarize  in Appendix C the 
fusions of the DVA current and the `$T$-system' for the 
algebra $A^{(1)}_{N-1}$.
\medskip

In the rest of this section, we briefly sketch 
the derivation of Proposition \ref{thm}.

\subsection{OPE's}
Let $r$ be generic, for a while.
We set 
\begin{eqnarray}
&&\rs=r-1, \\
&&\qi{u}={x^u-x^{-u}\over x-x^{-1}},\label{qi}\\
&& \ma{u}={\qi{u}\over \qi{ u+\rs+1}}={1\over \ma{-u-\rs-1}}.\label{ma}
\end{eqnarray}
In this section we prefer to use the additive notation 
and write $f(u)$ for the structure function $f(z)$ ($z=x^{2u}$)
in \eqref{fz}. 
It satisfies the relations:
\begin{lem}\label{str}
\begin{eqnarray*}
&(i)&{f(u-{1\over 2}) f(u+{1\over 2} )\over f(u)}=
{\ma{u-\rs-1/2}\over \ma{u-1/2}},\\
&(ii)& f(u)f(u\pm{3\over 2})=
{\ma{ \pm u- \rs}\over \ma{\pm u}}
{\ma{\pm u- \rs+ 1/2}\over \ma{\pm u+ 1/2}},\\
&(iii)&f(u-1)f(u) f(u+1 )
=
{\ma{u-\rs-1}\over \ma{u-1}}
{\ma{u-\rs-1/2}\over \ma{u-1/2}}
{\ma{u-\rs}\over \ma{u}}.
\end{eqnarray*}
\end{lem}
The operators $\Lambda_\pm(z),\Lambda_0(z)$ are defined by \eqref{Tz}.
\begin{lem}\label{ope}
The operator product expansions (OPE's) among $\Lambda_i(u)$ are 
\begin{eqnarray*}
&& f(u_2-u_1) \Lambda_i(x^{2u_1})\Lambda_j(x^{2u_2})\\
&=\!\!&
:\Lambda_i(x^{2u_1})\Lambda_j(x^{2u_2}): \\
&&\qquad\quad \times
\left\{
\begin{array}{ll}
  1&(i,j)=(+,+),\\
  {\dsp{\ma{u_1-u_2-\rs}\over\ma{u_1-u_2}}}&(i,j)=(+,0) ,\\
  {\dsp{\ma{u_1-u_2-\rs}\over\ma{u_1-u_2}}
  {\ma{u_1-u_2-\rs+1/2}\over\ma{u_1-u_2+1/2}}}&(i,j)=(+,-) ,\\
  {\dsp{\ma{u_1-u_2-\rs-1}\over\ma{u_1-u_2-1}}}&(i,j)=(0,+) ,\\
  {\dsp{\ma{u_1-u_2-\rs-1/2}\over\ma{u_1-u_2-1/2}}}&(i,j)=(0,0) ,\\
  {\dsp{\ma{u_1-u_2-\rs}\over\ma{u_1-u_2}}}&(i,j)=(0,-) ,\\
  {\dsp{\ma{u_1-u_2-\rs-1}\over\ma{u_1-u_2-1}}
  {\ma{u_1-u_2-\rs-3/2}\over\ma{u_1-u_2-3/2}}}&(i,j)=(-,+) ,\\
  {\dsp{\ma{u_1-u_2-\rs-1}\over\ma{u_1-u_2-1}}}&(i,j)=(-,0) ,\\
  1&(i,j)=(-,-). 
\end{array}\right. 
\end{eqnarray*}
\end{lem}

\subsection{Fused \dva currents $T_n(u)$ and $T_{\overline{n}}(u)$}
Suggested by the bootstrap program for general $r$,  
we introduce the following fused currents of \dvaa. 
\begin{dfn}\label{fused dva}
Define the fused DVA currents $T_n(u)$ and $T_{\overline{n}}(u)$ by
\begin{eqnarray*}
&&T_0(u)=T_{\overline{0}}(u)={\rm id},\quad 
T_1(u)=T_{\overline{1}}(u)=T(x^{2u}),\\
&&T_n(u)=\prod_{1\leq i<j\leq n} f(u_j-u_i)\cdot
T_1(u_1) T_1(u_2) \cdots T_1(u_n) 
\Biggl|_{u_i=u+{n-1\over 2}\rs -(i-1)\rs \atop
 \!\!\!\!\!\!\!\!\! \!\!\!\!\!\!\!\!\! \!\!\!\!\!\!\!\!\!
\!\!\!\!\!\!\!\!\!\!\!\!\!\!\!\!\!\! 1\leq i\leq n}, \\
&&T_{\overline{n}}(u)=
\prod_{1\leq i<j\leq n} f(u_j-u_i)\cdot
T_1(u_1) T_1(u_2) \cdots T_1(u_n) 
\Biggl|_{u_i=u+{n-1\over 2}(\rs-{1\over 2})  -(i-1)(\rs-{1\over 2}) \atop 
 \!\!\!\!\!\!\!\!\!\!\!\!\!\!\!\!\!\! \!\!\!\!\!\!\!\!\!
\!\!\!\!\!\!\!\!\!\!\!\!\!\!\!\!\!\! 
\!\!\!\!\!\!\!\!\!\!\!\!\!\!\!\!\!\! \!\!\!\!\!\!\! 
1\leq i\leq n}.
\end{eqnarray*}
Define 
the structure functions for $T_n(u)$, $T_{\overline{n}}(u)$ by
\begin{eqnarray*}
&&f_{mn}(u)=\prod_{i=1}^m\prod_{j=1}^n
f(u+{n-m\over 2}\rs -(j-i)\rs),\\
&& f_{\overline{m}\,\overline{n}}(u)=\prod_{i=1}^m\prod_{j=1}^n
f(u+{n-m\over 2}(\rs-{1\over 2}) -(j-i)(\rs-{1\over 2})),\\
&& f_{\overline{m}\,n}(u)=f_{n\,\overline{m}}(u)
=\prod_{i=1}^m\prod_{j=1}^n 
f(u+{n+1\over 2}\rs-
{m+1\over 2}(\rs-{1\over 2}) - j\rs+
i(\rs-{1\over 2})).
\end{eqnarray*}
\end{dfn}

These fused DVA currents enjoy the ZF exchange relations.
\begin{lem}\label{exchange}
As meromorphic functions, the following exchange relations hold.
\begin{eqnarray*}
&&f_{mn}(v-u)T_m(u)T_n(v)=f_{nm}(u-v)T_n(v)T_m(u),\\
&&f_{\overline{m}\,\overline{n}}(v-u)
T_{\overline{m}}(u)T_{\overline{n}}(v)=
f_{\overline{n}\,\overline{m}}(u-v)
T_{\overline{n}}(v)T_{\overline{m}}(u),\\
&&f_{m\,\overline{n}}(v-u)
T_{m}(u)T_{\overline{n}}(v)=
f_{\overline{n}\,m}(u-v)
T_{\overline{n}}(v)T_{m}(u).
\end{eqnarray*}
\end{lem}

\begin{lem}\label{anti}
\begin{eqnarray*}
&(i)&\ma{u-v \pm 1}f(v-u)T_1(u)T_1(v)\Biggr|_{u=v\mp 1}
= {\mp T_1(v\mp 1/2)\over \ma{-1/2}\ma{-1}} ,\\
&(ii)&\ma{u-v \pm {3/2}}f(v-u)T_1(u)T_1(v)\Biggr|_{u=v\mp 3/2}
= {\mp{\rm id}\over \ma{1/2}\ma{-1}\ma{-3/2}},\\
&(iii)&\ma{u_1-u_2-1}\ma{u_2-u_3-1}\\
&&\quad \times f(u_2-u_1)f(u_3-u_1)f(u_3-u_2)
T_1(u_1)T_1(u_2)T_1(u_3)\Biggr|_{u_2=u_1-1\atop u_3=u_2-1}\\
&&=
{{\rm id}\over 
\ma{1}\ma{{1\over2}}\ma{-{1\over2}}
\ma{-1}\ma{-1}\ma{-{3\over 2}}\ma{-2}}.
\end{eqnarray*}
\end{lem}

We need the analyticity properties of the operator products.
(Some details of the derivation 
are given in Appendix for the case of $A^{(1)}_{N-1}$.) 

\begin{lem}\label{lem1}
The product $f_{mn}(v-u)T_m(u)T_n(v)$ has poles only at 
\begin{eqnarray*}
u-v= 
\left\{
\begin{array}{l}
\left({m+n\over 2}-k\right)\rs -1, \\
\left({m+n\over 2}-k\right)\rs -{3\over 2}, \\
-\left({m+n\over 2}-k\right)\rs +1, \\
-\left({m+n\over 2}-k\right)\rs +{3\over 2}, \\
\end{array} \right.
\quad k=1,2,\cdots ,{\rm min}(m,n).
\end{eqnarray*}
All the poles are simple.
\end{lem}
\begin{lem}\label{lem2}
The product $f_{\overline{m}\,\overline{n}}(v-u)
T_{\overline{m}}(u)T_{\overline{n}}(v)$ has poles only at
\begin{eqnarray*}
&& \mbox{simple pole}:\\
&&u-v= 
\left\{
\begin{array}{l}
\left({m+n\over 2}-k\right)\left(\rs-{1\over 2}\right) -{3\over 2}, \\
-\left({m+n\over 2}-k\right)\left(\rs-{1\over 2}\right) +{3\over 2}, \\
\end{array} \right.
 k=1,2,\cdots ,{\rm min}(m,n), 
\end{eqnarray*}
and
\begin{eqnarray*}
&& \mbox{poles with multiplicity }
{\rm min}({\rm min}(m,n),{\rm min}(l,m+n-l)):\\
&&u-v= 
\left\{
\begin{array}{l}
\left({m+n\over 2}-l\right)\left(\rs-{1\over 2}\right) -1, \\
-\left({m+n\over 2}-l\right)\left(\rs-{1\over 2}\right) +1, \\
\end{array} \right.
 l=1,2,\cdots ,m+n-1.
\end{eqnarray*}
\end{lem}
\begin{lem}\label{lem31}
The product $f_{\overline{m}\,n}(v-u)T_{\overline{m}}(u)T_n(v)$ has
poles only at
\begin{eqnarray*}
u-v= 
\left\{
\begin{array}{l}
{n-1\over 2}\rs+
\left({m+1\over 2}-l\right)\left(\rs-{1\over 2}\right) -1, \\
{n-1\over 2}\rs+
\left({m+1\over 2}-k\right)\left(\rs-{1\over 2}\right) -{3\over 2}, \\
-{n-1\over 2}\rs
-\left({m+1\over 2}-l\right)\left(\rs-{1\over 2}\right) +1, \\
-{n-1\over 2}\rs
-\left({m+1\over 2}-k\right)\left(\rs-{1\over 2}\right) +{3\over 2}, \\
\end{array} \right.
\begin{array}{l}
 l=1,2,\cdots ,m, \\
 k=1,2,\cdots ,{\rm min}(m,n). 
\end{array}
\end{eqnarray*}
All the poles are simple.
\end{lem}

\subsection{The case of the dilute $A_3$ model}

The parameters for the dilute $A_3$ model are given by 
\begin{eqnarray*}
&&L=3,\qquad r=2{L+1\over L+2}={8 \over 5},\qquad \rs={3\over 5}.
\end{eqnarray*}

For $L=3$, we expect to have the following extra symmetry for $T_n(u)$.
\begin{conj}\label{con}
As operators acting on the BRST cohomology $H^{0}(C_{l,k})$ 
($k=1,2,3, l=1,2$), we have 
\begin{eqnarray*}
&(i)&
T_3(u) ={T_2(u) \over \ma{{2\over 10}}\ma{{3\over 10}}},\\
&(ii)&
T_4(u)=
{-T_1(u) \over
 \ma{-{9\over 10}}\ma{-{4\over 10}}
 \ma{-{3\over 10}}\ma{{3\over 10}}
 \ma{{2\over 10}}\ma{{8\over 10}} },\\
&(iii)&
T_5(u)=
{-{\rm id} \over 
 \ma{-{15\over 10}}
 \ma{-{10\over 10}}\ma{-{9\over 10}}
 \ma{-{4\over 10}}\ma{-{3\over 10}} \ma{{2\over 10}} 
 \ma{{3\over 10}}\ma{{8\over 10}}
 \ma{{9\over 10}}\ma{{14\over 10}} }.
\end{eqnarray*}
\end{conj}

One of the grounds for this conjecture \ref{con} is 
the degeneration of the structure functions.
\begin{lem}\label{lemstr}
For $\rs=3/5$, we have
\begin{eqnarray*}
&(i)&f_{13}(u)
=f_{12}(u)
{ \ma{u-{7\over 10}} \ma{u-{12\over 10}} \over
 \ma{u-{4\over 10}} \ma{u-{9\over 10}}} ,\\
&(ii)&f_{14}(u)=f_{11}(u)
{ \ma{u-{10\over 10}} \ma{u-{15\over 10}} \over
 \ma{u-{1\over 10}} \ma{u-{6\over 10}}} ,\\
&(iii)&f_{15}(u)=
{ \ma{u-{13\over 10}} \ma{u-{18\over 10}} \over
 \ma{u+{2\over 10}} \ma{u-{3\over 10}}} .
\end{eqnarray*}
\end{lem}

Using Lemma \ref{lemstr}, \ref{lem1}, \ref{lem2} and \ref{lem31}, 
we can check that 
the replacement $T_m(u) \leftrightarrow T_{5-m}(u)$ 
will not affect the analyticity in any of the OPE's acting on the
BRST cohomology space.

To obtain the correct proportionality constants 
in Conjecture \ref{con}, we calculated 
$\langle l,k| T_m(u) |l,k\rangle$ for $k=1,2,3, l=1,2$.

\subsection{Fusions of $T_2(u)$ at $\rs=3/5$}

If we study the bootstrap for $E_8$-symmetric particles 
carefully \cite{Zam89}, we realize that 
it is helpful to consider the fusions of $T_2(u)$. 

\begin{lem}\label{T68}
For $\rs=3/5$, the following equalities hold.
\begin{eqnarray*}
&(i)
&f_{22}(u_2-u_1)T_2(u_1)T_2(u_2)
\Biggr|_{u_1=u+{1\over 20} \atop
         u_2=u-{1\over 20}}\\
&& \qquad=
{\ma{-{8\over 10}}\ma{-{7\over 10}}\ma{-{2\over 10}}
     \ma{{3\over 10}}\over
\ma{-{5\over 10}}\ma{-{4\over 10}}\ma{-{3\over 10}}}
f_{1\overline{3}}(u_2-u_1)T_{1}(u_1)T_{\overline{3}}(u_2)
\Biggr|_{u_1=u+{11\over 20} \atop
         u_2=u-{3\over 20}}\\
&&\qquad=
{\ma{-{8\over 10}}\ma{-{7\over 10}}\ma{-{2\over 10}}
     \ma{{3\over 10}}\over
\ma{-{5\over 10}}\ma{-{4\over 10}}\ma{-{3\over 10}}}
f_{\overline{3}1}(u_2-u_1)T_{\overline{3}}(u_1)T_{1}(u_2)
\Biggr|_{u_1=u+{3\over 20} \atop
         u_2=u-{11\over 20}},\\
&(ii)
&f_{22}(u_2-u_1)T_2(u_1)T_2(u_2)
\Biggr|_{u_1=u+{7\over 20} \atop
         u_2=u-{7\over 20}} \\
&&\qquad=
{\ma{-{8\over 10}}\over
\ma{-{5\over 10}}\ma{-{3\over 10}}\ma{-{2\over 10}}}
f_{1\overline{2}}(u_2-u_1)T_{1}(u_1)T_{\overline{2}}(u_2)
\Biggr|_{u_1=u+{9\over 20} \atop
         u_2=u-{4\over 20}}\\
&&\qquad=
{\ma{-{8\over 10}}\over
\ma{-{5\over 10}}\ma{-{3\over 10}}\ma{-{2\over 10}}}
f_{\overline{2}1}(u_2-u_1)T_{\overline{2}}(u_1)T_{1}(u_2)
\Biggr|_{u_1=u+{4\over 20} \atop
         u_2=u-{9\over 20}}.
\end{eqnarray*}
\end{lem}

To prove these, we use 
\begin{eqnarray*}
&&{5\over 2}\rs -{3\over 2}=0,\qquad
T_2(u)=\ma{{2\over 10}}\ma{{3\over 10}}T_3(u),
\end{eqnarray*}
and Lemma \ref{anti}.

\subsection{$T_{\overline{5}}(u),
T_{\overline{6}}(u),T_{\overline{7}}(u)$ at $\rs=3/5$}

The fused DVA currents $T_{\overline{5}}(u),
T_{\overline{6}}(u),T_{\overline{7}}(u)$ for $\rs=3/5$ can be 
rewritten as follows.

\begin{lem}\label{kakikae} For $\rs=3/5$, we have
\begin{eqnarray*}
T_{\overline{5}}(u)
&\!\!=\!\!& \ma{-{1\over 2}}\ma{-1} \left(
{ \ma{-{6\over 10}}\ma{-{7\over 10}}\over 
\ma{-{12\over 10}}\ma{-{13\over 10}}} \right)^2
\ma{u_2'-u_2-1}\\
&&\times f_{\overline{2}1}(u_2-u_1) 
f_{\overline{2}1}(u_2'-u_1) 
f_{\overline{2}\,\overline{2}}(u_3-u_1) 
f_{11}(u_2'-u_2) f_{1\overline{2}}(u_3-u_2')
f_{1\overline{2}}(u_3-u_2)\\
&&\times T_{\overline{2}}(u_1)T_1(u_2)
T_1(u_2')T_{\overline{2}}(u_3)
\Biggr|_{u_1=u+{3\over 20}\atop
{{u_2=u-{1\over 2} \atop u_2'=u+{1\over 2} } \atop
u_3=u-{3\over 20}}} ,\\
T_{\overline{6}}(u)
&\!\!=\!\!&
\ma{-{1\over 2}}\ma{-1}
{ \ma{-{6\over 10}}\ma{-{7\over 10}} \over
\ma{-{12\over 10}}\ma{-{13\over 10}} } 
{\ma{-{6\over 10}}\ma{-{7\over 10}}
 \ma{-{8\over 10}} \over 
 \ma{-{12\over 10}}\ma{-{13\over 10}}
 \ma{-{14\over 10}} }
\ma{u_2'-u_2-1}\\
&&\times f_{\overline{3}1}(u_2-u_1) 
f_{\overline{3}1}(u_2'-u_1) 
 f_{\overline{3}\,\overline{2}}(u_3-u_1)
 f_{11}(u_2'-u_2)
 f_{1\overline{2}}(u_3-u_2')
f_{1\overline{2}}(u_3-u_2)\\
&&\times T_{\overline{3}}(u_1)T_1(u_2)
T_1(u_2')T_{\overline{2}}(u_3)
\Biggr|_{u_1=u+{3\over 20}\atop
{{u_2=u-{11\over 20} \atop u_2'=u+{9\over 20} } \atop
u_3=u-{4\over 20}}} ,\\
T_{\overline{7}}(u)
&\!\!=\!\!&
f_{\overline{5}2}(u_2-u_1) 
T_{\overline{5}}(u_1)T_2(u_2)\Biggr|_{u_1=u \atop u_2=u}\\
&\!\!=\!\!&
\ma{-{1\over 2}}\ma{-1}
\left(
{\ma{-{6\over 10}}\ma{-{7\over 10}}
 \ma{-{8\over 10}} \over 
 \ma{-{12\over 10}}\ma{-{13\over 10}}
 \ma{-{14\over 10}} } \right)^2
\ma{u_2'-u_2-1}\\
&&\times f_{\overline{3}1}(u_2-u_1) 
f_{\overline{3}1}(u_2'-u_1) 
 f_{\overline{3}\,\overline{3}}(u_3-u_1)
 f_{11}(u_2'-u_2) f_{1\overline{3}}(u_3-u_2')
f_{1\overline{3}}(u_3-u_2)\\
&&\times T_{\overline{3}}(u_1)T_1(u_2)
T_1(u_2')T_{\overline{3}}(u_3)
\Biggr|_{u_1=u+{4\over 20}\atop
{{u_2=u-{1\over 2} \atop u_2'=u+{1\over 2} } \atop
u_3=u-{4\over 20}}} .
\end{eqnarray*}
\end{lem}

\subsection{Structure functions}

In accordance with the currents $T^{(a)}(u)$, we introduce the
following structure functions.
\begin{dfn}\label{fg}
Define 
$f^{(a)}$ $(a=1,2,\cdots,8)$, $g^{(a)}$ $(a=2,3\cdots,6)$ by
\begin{eqnarray*}
&&f^{(a)}(u_1,u_2)=  f_{\overline{a}\,\overline{a}}(u_2-u_1)\qquad
(1\leq a\leq 5), \\
&&g^{(a)}(u_1,u_2)=  f_{\overline{a-1}\,\overline{a+1}}(u_2-u_1)\qquad
(2\leq a\leq 4), \\
&& g^{(5)}(u_1,u_2,u_3)= 
\ma{-{1\over 2}}\ma{-1}
{ \ma{-{6\over 10}}\ma{-{7\over 10}} \over
\ma{-{12\over 10}}\ma{-{13\over 10}} } 
{\ma{-{6\over 10}}\ma{-{7\over 10}}
 \ma{-{8\over 10}} \over 
 \ma{-{12\over 10}}\ma{-{13\over 10}}
 \ma{-{14\over 10}} }
{\ma{-{4\over 10}}\ma{-{5\over 10}}\ma{-{6\over 10}}
 \ma{-{7\over 10}} \over 
 \ma{-{10\over 10}} \ma{-{11\over 10}}\ma{-{12\over 10}}
 \ma{-{13\over 10}} } \\
&&\qquad\times
f_{\overline{4}\,\overline{3}}(u_2-u_1+{3\over 20})
f_{\overline{4}1}(u_2-u_1-{11\over 20}) 
f_{\overline{4}1}(u_3-u_1+{9\over 20})
f_{\overline{4}\,\overline{2}}(u_3-u_1-{4\over 20})\\
&&\qquad\times \ma{u_3-u_2}
f_{\overline{3}1}(u_3-u_2+{6\over 20})
f_{\overline{3}\,\overline{2}}(u_3-u_2-{7\over 20})\\
&&\qquad\times f_{11}(u_3-u_2+1)
f_{1\overline{2}}(u_3-u_2+{7\over 20}),\\
&& f^{(6)}(u_1,u_2)=
\ma{-{1\over 2}}\ma{-1}
\left(
{\ma{-{6\over 10}}\ma{-{7\over 10}}
 \ma{-{8\over 10}} \over 
 \ma{-{12\over 10}}\ma{-{13\over 10}}
 \ma{-{14\over 10}} } \right)^2 
{\qi{u_1-u_2-{1\over 10}}\over\qi{{16\over 10}}}
\ma{u_1-u_2-{1\over 10}}\\
&&\qquad\times f_{1\overline{3}}(u_2-u_1+{4\over 10}) 
f_{\overline{3}\,\overline{3}}(u_2-u_1-{3\over 10})
f_{11}(u_2-u_1+{11\over 10}) 
f_{\overline{3}1}(u_2-u_1+{4\over 10}), \\
&& g^{(6)}(u_1,u_2)= f_{\overline{5}2}(u_2-u_1) ,\\
&& f^{(7)}(u_1,u_2)=
{\ma{-{5\over 10}}\ma{-{4\over 10}}\ma{-{3\over 10}} \over 
\ma{-{8\over 10}}\ma{-{7\over 10}}\ma{-{2\over 10}}
     \ma{{3\over 10}}} \times f_{22}(u_2-u_1) ,\\
&& f^{(8)}(u_1,u_2)=\ma{-{1\over 2}}\ma{-1}
 \left(
{ \ma{-{6\over 10}}\ma{-{7\over 10}}\over 
\ma{-{12\over 10}}\ma{-{13\over 10}}} \right)^2
\ma{u_1-u_2-{1\over 10}}\\
&&\qquad\times 
f_{1\overline{2}}(u_2-u_1-{5\over 20}) 
f_{\overline{2}\,\overline{2}}(u_2-u_1+{8\over 20})
f_{11}(u_2-u_1-{18\over 20}) 
f_{\overline{2}1}(u_2-u_1-{5\over 20}).
\end{eqnarray*}
\end{dfn}

Collecting all the information together, 
we easily obtain the `$T$-system' with 
the $E_8$ symmetry stated in Proposition\ref{thm}.

\bigskip

\noindent
{\it Acknowledgments.}\quad 
We thank Atsuo Kuniba, Yaroslav Pugai and Junji Suzuki for 
discussions and interest. 
M. J. is grateful to Olivier Babelon and Jean Avan for 
kind invitation and hospitality during his stay in Paris VI, 
where a part of this work was done. 

\appendix
\setcounter{equation}{0}
\section{Operator product expansions}\lb{app:1}

We list the normal ordering relations. 
For operators $A(z),B(w)$ that have the form 
$:\exp(\mbox{linear in boson}):$, we use the notation
\be
A(z)B(w)=\dbr{A(z)B(w)}:A(z)B(w):
\en
and write down only the part $\dbr{A(z)B(w)}$. 
\begin{eqnarray*}
&&
\dbr{x_+(z_1)x_+(z_2)}=z_1^{\frac{r}{r-1}}(1-z_2/z_1)
\frac{(x^{-2}z_2/z_1,x^{2r-1}z_2/z_1;x^{2r-2})_\infty}
{(x^{-1}z_2/z_1,x^{2r}z_2/z_1;x^{2r-2})_\infty},
\\
&&\dbr{x_-(z_1)x_-(z_2)}=z_1^{\frac{r-1}{r}}(1-z_2/z_1)
\frac{(x^{2}z_2/z_1,x^{2r-1}z_2/z_1;x^{2r})_\infty}
{(xz_2/z_1,x^{2r-2}z_2/z_1;x^{2r})_\infty},
\\
&&
\dbr{x_\pm(z_1)x_\mp(z_2)}
=z_1^{-1}
\frac{1+z_2/z_1}
{(1+xz_2/z_1)(1+x^{-1}z_2/z_1)},
\\
&&
\dbr{\Phi_-(z_1)x_-(z_2)}
=z_1^{-\frac{r-1}{r}}
\frac{(x^{2r-1}z_2/z_1;x^{2r})_\infty}{(xz_2/z_1;x^{2r})_\infty},
\\
&&
\dbr{x_-(z_2)\Phi_-(z_1)}
=z_2^{-\frac{r-1}{r}}
\frac{(x^{2r-1}z_1/z_2;x^{2r})_\infty}{(xz_1/z_2;x^{2r})_\infty},
\\
&&
\dbr{\Phi_-(z_1)x_+(z_2)}
=x_+(z_2)\Phi_-(z_1)
=(z_1+z_2),
\\
&&
\dbr{\Psi^*_-(z_1)x_+(z_2)}
=z_1^{-\frac{r}{r-1}}
\frac{(x^{2r-1}z_2/z_1;x^{2r-2})_\infty}{(x^{-1}z_2/z_1;x^{2r-2})_\infty}
,
\\
&&
\dbr{x_+(z_2)\Psi_-^*(z_1)}
=z_2^{-\frac{r}{r-1}}
\frac{(x^{2r-1}z_1/z_2;x^{2r-2})_\infty}
{(x^{-1}z_1/z_2;x^{2r-2})_\infty},
\\
&&
\dbr{\Psi^*_-(z_1)x_-(z_2)}
=\dbr{x_-(z_2)\Psi^*_-(z_1)}
=(z_1+z_2),
\\
&&\dbr{\Phi_-(z_1)\Phi_-(z_2)}
=z_1^{\frac{r-1}{r}}
\frac{(x^2z_2/z_1,x^3z_2/z_1,
x^{2r+3}z_2/z_1,x^{2r+4}z_2/z_1;x^6,x^{2r})_\infty}
{(x^5z_2/z_1,x^6z_2/z_1,
x^{2r}z_2/z_1,x^{2r+1}z_2/z_1;x^6,x^{2r})_\infty},
\\
&&
\dbr{\Psi^*_-(z_1)\Psi^*_-(z_2)}
=z_1^{\frac{r}{r-1}}
\frac{(z_2/z_1,xz_2/z_1,
x^{2r+3}z_2/z_1,x^{2r+4}z_2/z_1;x^6,x^{2r-2})_\infty}
{(x^3z_2/z_1,x^4z_2/z_1,
x^{2r}z_2/z_1,x^{2r+1}z_2/z_1;x^6,x^{2r-2})_\infty},
\\
&&
\dbr{\Phi_-(z_1)\Psi^*_-(z_2)}
=z_1^{-1}
\frac{(-x^4z_2/z_1,-x^5z_2/z_1;x^6)_\infty}
{(-xz_2/z_1,-x^2z_2/z_1;x^6)_\infty},
\\
&&
\dbr{\Psi^*_-(z_2)\Phi_-(z_1)}
=z_2^{-1}
\frac{(-x^4z_1/z_2,-x^5z_1/z_2;x^6)_\infty}
{(-xz_1/z_2,-x^2z_1/z_2;x^6)_\infty}.
\end{eqnarray*}

As meromorphic functions we have
\begin{eqnarray*}
&&
x_+(z_1)x_+(z_2)=
\frac{[u_1-u_2+1]^*}{[u_1-u_2-1]^*}
\frac{[u_1-u_2-1/2]^*}{[-u_1+u_2-1/2]^*}
x_+(z_2)x_+(z_1),
\\
&&x_-(z_1)x_-(z_2)=
\frac{[u_1-u_2-1]}{[u_1-u_2+1]}
\frac{[u_1-u_2+1/2]}{[-u_1+u_2+1/2]}
x_-(z_2)x_-(z_1),
\\
&&x_\pm(z_1)x_\mp(z_2)=
x_\mp(z_2)x_\pm(z_1),
\\
&&\Phi_-(z_1)x_-(z_2)=
\frac{[u_1-u_2+1/2]}{[-u_1+u_2+1/2]}
x_-(z_2)\Phi_-(z_1),
\\
&&
\Phi_-(z_1)x_+(z_2)=x_+(z_2)\Phi_-(z_1),
\\
&&
\Psi^*_-(z_1)x_+(z_2)=
\frac{[u_1-u_2-1/2]^*}{[-u_1+u_2-1/2]^*}
x_+(z_2)\Psi^*_-(z_1),
\\
&&
\Psi^*_-(z_1)x_-(z_2)=x_-(z_2)\Psi^*_-(z_1),
\\
&&\Phi_-(z_1)\Phi_-(z_2)=
\rho(u_2-u_1)\Phi_-(z_2)\Phi_-(z_1),
\\
&&
\Psi^*_-(z_1)\Psi^*_-(z_2)=
\rho^*(u_1-u_2)\Psi^*_-(z_2)\Psi^*_-(z_1),
\\
&&\Phi_-(z_1)\Psi^*_-(z_2)=
\tau(u_2-u_1)
\Psi^*_-(z_2)\Phi_-(z_1).
\end{eqnarray*}
Here $\rho(u)$, $\rho^*(u)$ and $\tau(u)$ 
are given respectively by \eqref{rho},\eqref{rho*} and \eqref{tau}.

\setcounter{equation}{0}
\section{BRST charges}\lb{app:2}

We give here a proof of the properties of BRST charges
stated in subsection \ref{subsec:3.5}.  
The method is a proper adaptation of 
the work \cite{FJMOP} to the present situation. 

\subsection{Feigin-Odesskii algebra}

First let us prepare the notation. 
Let $\At_n$ be the set of all functions $F(u_1,\cdots,u_n)$
which is holomorphic on $\C^n$, symmetric in $u_1,\cdots,u_n$, and enjoys
the quasi-periodicity properties ($\rs=r-1$) 
\begin{eqnarray}
&&F(u_1+\rs,u_2,\cdots,u_n)=(-1)^n F(u_1,u_2,\cdots,u_n),
\lb{quasi1}\\
&&F(u_1+\tau,u_2,\cdots,u_n)=(-1)^n F(u_1,u_2,\cdots,u_n)\nn\\
&&\qquad\qquad\qquad\qquad\qquad\qquad\times 
\exp\Bigl(\frac{2\pi i}{\rs}(n u_1-\sum_{j=2}^n u_j 
-\frac{n-1}{2}+n\frac{\tau}{2})\Bigr), 
\lb{quasi2}
\end{eqnarray}
where $\tau=\pi i/\log x$. 
Clearly 
\[
\At_1=\C f_1, \quad f_1(u)=[u]^*. 
\]
We have also $\dim \At_2=2$. 
If $F\in\At_n$ is not identically $0$, then 
it has $n$ zeroes $\{u_1^{(j)}\}_{j=1}^n$ $\bmod~\Z\rs\oplus\Z\tau$ 
satisfying 
\[
\sum_{j=1}^n u_1^{(j)}=\sum_{j=2}^n u_2+\frac{n-1}{2}.
\]

Let $F\in \At_m$, $G\in \At_n$. 
Following the line of \cite{FJMOP}, 
we define the $*$-product $F*G\in\At_{m+n}$ by 
\begin{eqnarray*}
(F*G)(u_1,\cdots,u_{m+n})
&\!\!=\!\!&\mbox{Sym}\Bigl(F(u_1-n,\cdots,u_m-n)G(u_{m+1},\cdots,u_{m+n})
\lb{*prod}\\
&&\qquad\quad \times\!\!
\prod_{1\le i\le m\atop m+1\le j\le m+n}
\frac{[u_i-u_j+1]^*[u_i-u_j-1/2]^*}{[u_i-u_j]^*}
\Bigr).
\nn
\end{eqnarray*}
Here the symbol $Sym$ stands for the symmetrization.
$\At=\oplus_{n=0}^\infty \At_n$ equipped with the $*$-product 
is an associative graded algebra with unit. 
We denote by $A=\oplus_{n=0}^\infty A_n$ ($A_n=\At_n\cap A$) 
the subalgebra of $\At$ generated by $\At_1$ and $\At_2$.

Let us say that a function $F(u_1,\cdots,u_n)$ has the 
property (P) if either $n=1,2$, or $n\ge 3$ and 
\[
F(u_1,\cdots,u_n)\Bigl|_{u_j-u_i=u_k-u_j=1/2}
\equiv 0
\quad \mbox{for any distinct $i,j,k$}.
\]
\begin{lem}\lb{lem3}
\begin{enumerate}
\item Elements of $A$ has the property (P). 
\item $A$ is commutative. 
\end{enumerate}
\end{lem}
\medskip

\noindent Proof.\quad 
{}From the definition of $*$ we can verify that,  
if $F\in\At_m$ and $G\in\At_n$ have the property (P), 
then so does $F*G$. Hence (i) follows.

Let $g\in A_2$ be an element linearly independent from $f_1*f_1\in A_2$.   
To see (ii), it suffices to show that $f_1*g=g*f_1$. 
Set $h=f_1*g-g*f_1$.  
A simple check shows that $h(u+1,u,u-1)=0$. 
By symmetry and the property (P), $h(u_1,u_2,u_2-1)$ 
viewed as a function of $u_1$ has zeroes at $u_2+1,u_2-1/2,u_2-2$. 
Since their sum is different from $2u_2\bmod~\Z\rs\oplus\Z\tau$, 
we have $h(u_1,u_2,u_2-1)\equiv 0$. 
By symmetry, this implies that $h(u_1,u_2,u_3)$ 
has $4$ zeroes $u_1=u_2\pm 1,u_3\pm 1$.  Hence $h\equiv 0$. 
\qed

\subsection{BRST charges}

For $F\in A_n$, we define 
\begin{eqnarray*}
Q(F)&\!\!=\!\!&\oint\!\cdots\!\oint 
\prod_{j=1}^n\frac{\dz_j}{[u_j+\frac{1}{2}]^*}~x_+(z_1)\cdots x_+(z_n) \\
&&\times 
\left(\prod_{1\le i<j\le n}\frac{[u_i-u_j]^*}{[u_i-u_j+1]^*[u_i-u_j-1/2]^*}
\right)
F(u_1+\hat{l}/2,\cdots,u_n+\hat{l}/2),
\end{eqnarray*}
where the contour is $|z_1|=\cdots=|z_n|=1$. 
Because of the quasi-periodicity \eqref{quasi2},
the integrand is single valued. 
Using the exchange relation 
\[
x_+(z_1)x_+(z_2)=-\frac{[u_1-u_2+1]^*}{[u_2-u_1+1]^*}
\frac{[u_1-u_2-1/2]^*}{[u_2-u_1-1/2]^*}
x_+(z_2)x_+(z_1)
\]
along with $\hat{l}\,x_+(z)=x_+(z)(\hat{l}-2)$, 
we find 
\begin{equation}
Q(F)Q(G)=Q(F*G).
\lb{FG}
\end{equation}

Now set
\begin{eqnarray*}
f_1(u)&\!\!=\!\!&[u]^*,\\
f^{(a)}_2(u_1,u_2)&\!\!=\!\!&[2a+1]^*[a-1/2]^*[u_1-a]^*[u_2+a-1]^*
[u_1-u_2+a-1/2]^*
\\
&&\!\!\!\!-[2a-1]^*[a+1/2]^*[u_1+a]^*[u_2-a-1]^*[u_1-u_2-a-1/2]^*, 
\end{eqnarray*}
and introduce $h_l\in A_l$ ($1\le l\le L$) by  
\begin{eqnarray*}
&&h_{2m+1}=f_1*f^{(1)}_2*\cdots*f^{(m)}_2,
\qquad \Bigl(0\le m\le \frac{L-1}{2}\Bigr),\\
&&h_{2m}=f^{((L+1)/2-m)}_2*\cdots*f^{((L-3)/2)}_2*f^{((L-1)/2)}_2,
\qquad \Bigl(1\le m\le \frac{L-3}{2}\Bigr).
\end{eqnarray*}
We have $h_L=h_{l}*h_{L-l}$ ($1\le l\le L-1$).
We define the BRST charges by 
\[
Q_l=Q(h_l). 
\]

Propositions \ref{prop:3.51},\ref{prop:3.52} reduce to the following 
assertions.
\begin{prop}\lb{prop:app3.1} 
$h_l(u_1,\cdots,u_l)$ can be written as 
\begin{eqnarray}
&&h_l(u_1,\cdots,u_l)=\bar{h}_l(u_1,\cdots,u_l)\times
\cases{{\dsp\prod_{i=1}^l\Bigl[u_i+\frac{-l+1}{2}\Bigr]^*} 
& (l : \mbox{\it odd}) \cr 
{\dsp\prod_{i=1}^l\Bigl[u_i+\frac{L-l+1}{2}\Bigr]^*} 
& (l : \mbox{\it even}), \cr }
\lb{h01}
\end{eqnarray}
where $\bar{h}_l$ are holomorphic and satisfies
\be
&&\bar{h}_l(u_1+\rs,\cdots,u_l)=(-1)^{l-1}\bar{h}_l(u_1,\cdots,u_l),
\\
&&\bar{h}_l(u_1+\tau,\cdots,u_l)=\bar{h}_l(u_1,\cdots,u_l)
\times\exp\Bigl(\frac{2\pi i}{\rs}\sum_{j=2}^l(u_1-u_j+\frac{\tau}{2})\Bigr).
\en
Moreover it is translationally invariant, i.e., 
\begin{eqnarray*}
&&\bar{h}_l(u_1+v,\cdots,u_l+v)=\bar{h}_l(u_1,\cdots,u_l).
\lb{h02}
\end{eqnarray*}
\end{prop}

\begin{prop}\lb{prop:app3.2} 
We have $h_L\equiv 0$. 
\end{prop}

\noindent Proof of Proposition \ref{prop:app3.2}.\quad 
Let $0\le m\le \frac{L-3}{2}$. 
In the equality $h_L=h_{2m+1}*h_{L-2m-1}$ we set $u_1=-m-1$. 
Using 
Proposition \ref{prop:app3.1}, 
we find that each summand in the symmetrization \eqref{*prod} vanishes. 
Similarly, if we set $u_1=m$, then each summand of 
$h_L=h_{L-2m-1}*h_{2m+1}$ vanishes. 
Therefore $h_L$ has $L-1$ zeroes 
$u_1=0,1,\cdots,\frac{L-3}{2},-1,-2,\cdots,-\frac{L-1}{2}$.

Suppose $h_L$ did not vanish identically. 
{}From the quasi-periodicity, $h_L$ has a zero at 
$u_1=\sum_{j=2}^Lu_j+L-1$. 
By symmetry, $u_1=u_2-\sum_{j=3}^Lu_j-(L-1)$
must also be a zero. This is a contradiction. 
\qed

In the next subsection we prove Proposition \ref{prop:app3.1}. 

\subsection{Proof of Proposition \ref{prop:app3.1}}

We prove Proposition \ref{prop:app3.1} for odd $l=2m+1$. 
The statement is obvious for $m=0$. 
Assuming $m\ge 1$ we proceed by induction on $m$.

\begin{lem}\lb{lem4}
For $m=1,\cdots,\frac{L-1}{2}$ we have 
\begin{eqnarray}
&&h_{2m+1}(m,m\pm 1, u_3,\cdots,u_{2m+1})\equiv 0. 
\lb{eq1}
\end{eqnarray}
\end{lem}
\medskip

\noindent Proof. \quad 
We use the property
\begin{equation}
f^{(a)}_2(\pm a,\pm a+1)=0.
\lb{fa}
\end{equation}
In the definition of $h_{2m+1}=h_{2m-1}*f^{(m)}_2$, 
set $u_1=m,u_2=m+1$. 
Using the induction hypothesis 
$h_{2m-1}(m-1,\cdots)\equiv 0$ and $f^{(m)}_2(m,m+1)=0$, 
we see that each summand vanishes.
Similarly if we set  $u_1=m,u_2=m-1$ 
in $h_{2m+1}=f^{(m)}_2*h_{2m-1}$ and  
use  $f^{(m)}_2(-m,-m+1)=0$, the result is zero. 
\qed

\begin{lem}\lb{lem5}For $t=2,3,\cdots$, we have
\begin{equation}
h_{2m+1}(m,u_2,\cdots,u_{2m+1})
\Bigl|_{u_{2s+1}=u_{2s}+1/2 \atop t\le s\le m}
\equiv 0.
\lb{h2}
\end{equation}
\end{lem}
Taking $t=m+1$ we obtain
the first assertion of Proposition \ref{prop:app3.1}.
\medskip

\noindent Proof. \quad 
Denote the left hand side of \eqref{h2} by $g_t$. 
We show $g_t\equiv 0$ by induction on $t$. 

Let $t=2$, and consider first 
\begin{equation}
h_{2m+1}(m,m\pm 1/2,u_3,\cdots,u_{2m+1})
\Bigl|_{u_{2s+1}=u_{2s}+1/2 \atop 2\le s\le m}.
\lb{eq2}
\end{equation}
As a function of $u_3$, \eqref{eq2} has $2m+1$ zeroes at 
$m+1,m-1, m\mp 1/2$ and $u_{2s}+1,u_{2s}-1/2$ ($2\le s\le m$). 
Comparing with the quasi-periodicity, we conclude that 
\eqref{eq2} vanishes identically. 
This means that $g_2$ as a function of $u_2$ has 
$2m+2$ zeroes at 
$m\pm 1,m\pm1/2$ and $u_{2s}+1,u_{2s}-1/2$ ($2\le s\le m$). 
Therefore $g_2\equiv 0$. 

Suppose we have shown $g_t\equiv 0$, and consider 
\begin{equation}
g_{t+1}=h_{2m+1}(m,u_2,\cdots, u_{2t},u_{2t+1},
u_{2t+2},u_{2t+2}+1/2,\cdots,u_{2m},u_{2m}+1/2).
\lb{h3}
\end{equation}
{}From the induction hypothesis, it vanishes for $u_{2t+1}=u_{2t}+1/2$. 
By symmetry it vanishes for 
\[
u_2=u_3\pm1/2,\cdots,u_{2t+1}\pm 1/2.
\]
It also vanishes for $u_2=m\pm 1$ and $u_2=u_{2s}-1/2,u_{2s}+1$ 
($t+1\le s\le m$). 
Since the number of zeroes exceed $2m+1$, 
we conclude $g_{t+1}\equiv 0$.
\qed

In the case of even $l$, we note that 
\[
f^{((L-1)/2)}_2(u_1,u_2)=[2]^*\Bigl[\frac{L}{2}\Bigr]^*
\Bigl[u_1+\frac{L-1}{2}\Bigr]^*\Bigl[u_2+\frac{L-1}{2}\Bigr]^*
\Bigl[u_1-u_2-\frac{L}{2}\Bigr]^*  
\]
has a zero at $u_1=(1-L)/2$. 
Using this the proof goes similarly. 

\begin{lem}\lb{lem6}
The function $\bar{h}_l$ in \eqref{h01} is translationally invariant. 
\end{lem}
\medskip

\noindent Proof.\quad
Consider
\[
\bar{h}_l(u_1+v,\cdots,u_l+v).
\]
It is holomorphic and doubly periodic in $v$, hence is a constant. 
The conclusion follows by setting $v=0$. 
\qed

\setcounter{equation}{0}
\section{Deformed $W$ algebra for $\slNBig$}\lb{app:3}

We discuss here the fusion of the deformed $W$ algebra (DWA) 
associated with $\slN$ \cite{qWN,FeFr95}. 

\subsection{Basic current}

Fix complex numbers $x,r^*\in\C$, $0<|x|<1$. 
We keep the notation $[u]_x$ \eqref{qi} and $\ma{u}$ \eqref{ma}.
In this appendix we consider the case of `generic' $\rs$, 
i.e., we assume that $m,n\in\Z$, $[m+nr^*]_x=0$ implies $m=n=0$. 

In the free field realization, 	
the simplest current of DWA for the algebra $\slN$ is presented in the form
\be
W_{(1)}(u)=\sum_{i=1}^N\Lambda_i(u). 
\en
Each $\Lambda_i(u)$ is a normally-ordered exponential of bosonic oscillators. 
Their explicit formula is irrelevant here
(see e.g. \cite{qWN}, eq.(2), wherein 
$z=x^{2u}$, $q=x^{2r^*+2}$, $t=x^{2r^*}$ in the present notation). 
We need only the following 
normal ordering rule for their products: 
\bea
f(u,v)\Lambda_i(u)\Lambda_j(v)=\,:\Lambda_i(u)\Lambda_j(v):\times
\left\{
\begin{array}{ll}
{\dsp{\ma{u-v-1-\rs}\over \ma{u-v-1}}} & (i<j) \\
1 & (i=j)\;, \\
{\dsp{\ma{u-v-\rs}\over \ma{u-v}}} & (i>j) 
\end{array} \right.
\lb{contra}
\ena
where the structure function $f(u,v)=f(v-u)$ is given by
$$
f(u)=\frac{1}{(1-x^{2u})}
\frac{(x^{2(u+N-1)},x^{2(u+1+r^*)},x^{2(u-r^*)};x^{2N})_\infty}
{(x^{2(u+1)},x^{2(u+N+r^*)},x^{2(u+N-r^*-1)};x^{2N})_\infty}.
$$

\begin{lem}\lb{lem:app21}
We have the exchange relation as meromorphic functions 
\bea
f(u,v)W_{(1)}(u)W_{(1)}(v)
=
f(v,u)W_{(1)}(v)W_{(1)}(u).
\lb{fww}
\ena
Both sides are regular except for simple poles 
at $u-v=\pm 1 \bmod \Gamma$, where $\Gamma=(\pi i/\log x)\Z$. 
\end{lem}
Notice that the pole $u=v$ which appears in 
\eqref{contra} is canceled in \eqref{fww}.
In general, each matrix element of the product 
\[
\prod_{1\le s<t\le m}f(u_s,u_t)\times
W_{(1)}(u_1)\cdots W_{(1)}(u_m)
\]
is a rational function of $x^{2u_i}$ with at most simple poles 
at $u_i-u_j=\pm 1$ ($i<j$). 
 
\subsection{Fused currents}

Let $\lambda=(\lambda_1,\cdots,\lambda_l)$ 
($\lambda_1\ge\cdots\ge\lambda_l>0$) be a partition. 
We identify $\lambda$ with a Young diagram.  
For $j,s=1,2,\cdots$, we attach a variable $u(j,s)$ 
to the box on the $j$-th row and $s$-th column of $\lambda$:
$$
\setlength{\unitlength}{1mm}
\begin{picture}(140,35)(-15,5)
\put(-5,15){\makebox(10,10)[l]{$\lambda\;=$}}
\put(10,8){\line(0,1){24}}
\put(10,32){\line(1,0){33}}
\put(10,8){\line(1,0){17}}
\put(27,8){\line(0,1){4}}
\put(27,12){\line(1,0){4}}
\put(31,12){\line(0,1){4}}
\put(31,16){\line(1,0){4}}
\put(35,16){\line(0,1){4}}
\put(35,20){\line(1,0){4}}
\put(39,20){\line(0,1){4}}
\put(39,24){\line(1,0){4}}
\put(43,24){\line(0,1){8}}
\put(20,18.5){\line(1,0){3}}
\put(20,18.5){\line(0,1){3}}
\put(20,21.5){\line(1,0){3}}
\put(23,18.5){\line(0,1){3}}
\put(5,17.5){\makebox(5,5){$j$}}
\put(19,32){\makebox(5,5){$s$}}
\multiput(11,20)(2,0){5}{\line(1,0){1}}
\multiput(21.5,31)(0,-2){5}{\line(0,-1){1}}
\put(45,17.5){\makebox(5,5)[b]{,}}
\put(55,15){\makebox(40,10)[l]{$u(j,s)=u-(j-1)-(s-1)\rs$.}}
\end{picture}
$$
For partitions $\la=(\la_1,\cdots,\la_l)$, 
$\mu=(\mu_1,\cdots,\mu_m)$, we set 
\bea
f_{\la,\mu}(u,v)=\prod_{1\le j\le l\atop 1\le s\le \la_j} 
\prod_{1\le k\le m\atop 1\le t\le \mu_k} 
f(u(j,s),v(k,t)).  
\lb{flamu}
\ena

We shall associate `fused' currents $W_\lambda(u)$ 
with each $\lambda$.  
First consider the case of a single row diagram $\la=(m)$.
\begin{dfn}\lb{dfn:app21}
\be
&&W_{(m)}(u)=
\Bigl(
\prod_{1\le s<t\le m}f(u_s,u_t)\times
W_{(1)}(u_1)\cdots W_{(1)}(u_m)
\Bigr)
\Bigl|_{u_s=u-(s-1)\rs \atop 1\le s\le m}.
\en
\end{dfn}
In view of the remark after Lemma \ref{lem:app21}, 
the right hand side is well-defined. 
Alternatively $W_{(m)}(u)$ can be defined inductively as 
\bea
W_{(m)}(u)&\!\!=\!\!&f_{(m-1),(1)}(u,u')W_{(m-1)}(u)W_{(1)}(u')
\Bigl|_{u'=u-(m-1)\rs}
\lb{sym1}\\
&\!\!=\!\!&f_{(1),(m-1)}(u,u')W_{(1)}(u)W_{(m-1)}(u')
\Bigl|_{u'=u-\rs}.
\lb{sym2}
\ena

\begin{lem}\lb{lem:app22}
We have 
\bea
&&f_{(1),(m)}(u,v)W_{(1)}(u)W_{(m)}(v)
=f_{(m),(1)}(v,u)W_{(m)}(v)W_{(1)}(u).
\lb{fww1}
\ena
Both sides of \eqref{fww1} are regular except 
for simple poles ($\bmod \Gamma$) at $u-v=-1, 1-(m-1)\rs$.
\end{lem}

\proof
The exchange relation \eqref{fww1} is obvious. 
Let us verify the statement about the position of poles
by induction on $m$. 
The case $m=1,2$ can be verified by direct calculation. 
Suppose it is true for $m-1$. 
Using the expression \eqref{sym1} and \eqref{sym2} 
and the induction hypothesis, we see that the possible poles 
in $u$ are confined to 
\be
\!\!\!\!\!
&&\{v-1,v+1-(m-2)\rs,v\pm1-(m-1)\rs\}\cap\{v\pm 1,v-\rs-1,v+1-(m-1)\rs\}
\\
\!\!\!\!\!&&=\{v-1,v+1-(m-1)\rs\}.
\en
\qed

Arguing similarly, we have
\begin{lem}\lb{lem:app221}
\bea
&&f_{(m),(n)}(u,v)W_{(m)}(u)W_{(n)}(v)
=f_{(n),(m)}(v,u)W_{(n)}(v)W_{(m)}(u).
\lb{fww11}
\ena
Both sides of \eqref{fww11} are regular except 
for simple poles ($\bmod \Gamma$) at 
\be
u-v&\!\!=\!\!&1-j\rs \qquad (\max(0,n-m)\le j\le n-1), \\
&\!\!=\!\!&-1+j\rs \qquad (\max(0,m-n)\le j\le m-1). 
\en
\end{lem}

For a general partition $\lambda=(\lambda_1,\cdots,\lambda_l)$,  
we define
\begin{dfn}\lb{dfn:app22}
\be
W_{\la}(u)=\Bigl(\prod_{i=1}^{l-1}\ma{u_i-u_{i+1}-1}\,\cdot\!\!\!
\prod_{1\le i<j\le l}\!\!f_{(\la_i),(\la_j)}(u_i,u_j)\cdot
W_{(\la_1)}(u_1)\cdots W_{(\la_l)}(u_l)
\Bigr)\Bigl|_{u_i=u-(i-1)\atop 1\le i\le l}.
\en
\end{dfn}
This definition makes sense by Lemma \ref{lem:app22}. 
We have 
\bea
f_{\la,\mu}(u,v)W_{\la}(u)W_{\mu}(v)
=W_{\mu}(v)W_{\la}(u)f_{\mu,\la}(v,u).
\lb{fww2}
\ena

In the case of a single column diagram $\la=(1^a)$, 
 $W_{(1^a)}(u)$ coincides with 
the fundamental DWA currents $W_a(z)$ in \cite{qWN,FeFr95} 
up to a numerical factor and a shift of $u$ (see \eqref{(1m)} below).  

We remark that another fused currents can be constructed similarly 
by replacing $u(j,s)=u-(j-1)-(s-1)r^*$ with $u-(j-1)+(s-1)r$ ($r=r^*+1$).

\subsection{Tableaux sum}

Let $\la$ be a partition. 
Denote by $SST(\la)$ the set of 
semi-standard tableaux of shape $\la$
on the letters $\{1,2,\cdots,N\}$. 
For $T\in SST(\la)$, we set 
\be
\Lambda_T(u)=\,:\!\!\!\prod_{1\le j\le l\atop 1\le s\le \la_j}
\Lambda_{T(j,s)}\Big(u(j,s)\Bigr):, 
\en
where $T(j,s)\in\{1,\cdots,N\}$ signifies 
the letter in the $(j,s)$-th position of $T$. 

The current $W_{\la}(u)$ is given explicitly as follows.
\begin{prop}
We have
\be
&&W_\la(u)=d_\la\sum_{T\in SST(\la)}c_T\cdot \Lambda_T(u).
\en
The coefficients $d_\la,c_T$ are given by
\be
d_\la&\!\!=\!\!&
\prod_{1+j<k}\frac{\ma{k-j-1-\la_j\rs}_{\la_k}}{\ma{k-j-1}_{\la_k}}\cdot
\prod_{k=2}^{l}\frac{\ma{-\la_{k-1}\rs}_{\la_{k}}}{\ma{\rs}_{\la_{k}-1}},
\\
c_{T}&\!\!=\!\!&
\prod_{j=1}^l\frac{\prod_{i=1}^N\ma{-1}_{w_{ji}}}{\ma{-1}_{\la_j}}\cdot
\prod_{j<k}\frac{\ma{k-j-\la_j\rs}_{\la_k}}
{\ma{k-j-1-\la_j\rs}_{\la_k}}
\\
&&\quad\times
\prod_{j<k}\prod_{i=1}^N
\frac{\ma{k-j-1+(s_{k,i-1}-s_{j,i-1})\rs}_{w_{k,i}}}
{\ma{k-j+(s_{k,i-1}-s_{j,i})\rs}_{w_{k,i}}}, 
\en
where 
$w_{ji}$ is the number of the letter $i$ in the $j$-th row of $T$ 
($1\le j\le l, 1\le i\le N$), $s_{ji}=w_{j1}+\cdots+w_{ji}$, and  
\be
\ma{u}_n=\ma{u}\ma{u+\rs}\cdots\ma{u+(n-1)\rs}.
\en
\end{prop}
We omit the proof. 
Notice that $c_{T_0}=1$ 
for the tableau $T_0$ with $T_0(j,s)=j$ for all $j,s$. 

\example 
\bea
W_{(m)}(u)=\sum_{w_1,\cdots,w_N\ge 0\atop w_1+\cdots+w_N=m}
\frac{\prod_{i=1}^N\ma{-1}_{w_{i}}}{\ma{-1}_{m}}
\cdot \Lambda_T(u),
\lb{(m)}
\ena
where $T=(1^{w_1},2^{w_2},\cdots,N^{w_N})$. 
\bea
&&W_{(1^a)}(u)=d_{(1^a)}\sum_{1\le i_1<\cdots<i_a\le N}
:\Lambda_{i_1}(u)\cdots\Lambda_{i_a}(u-a+1):.
\lb{(1m)}
\ena

\subsection{$W_\lambda(u)$ in terms of $W_{(1^a)}(u)$}

$W_{\la}(u)$ can also be obtained from $W_{(1^a)}(u)$. 
First note the following fact which 
can be shown similarly as Lemma \ref{lem:app221}.
\begin{lem}\lb{lem:app23}
We have 
\bea
&&\tilde{f}_{(1^a),(1^b)}(u,v)W_{(1^a)}(u)W_{(1^b)}(v)
=\tilde{f}_{(1^b),(1^a)}(v,u)W_{(1^b)}(v)W_{(1^a)}(u),
\lb{fww3}
\ena
where  $\tilde{f}_{\la,\mu}(u,v)$ is defined similarly as in \eqref{flamu} 
with $f(u,v)$ replaced by 
\be
\tilde{f}(u,v)=\frac{\ma{u-v-1}}{\ma{u-v-\rs}}f(u,v)
=\frac{\ma{v-u-1}}{\ma{v-u-\rs}}f(u,v).
\en 
Both sides of \eqref{fww3} 
are regular except for simple poles  ($\bmod \Gamma$) at 
\be
u-v&\!\!=\!\!&\rs-j \qquad (\max(0,b-a)\le j\le b-1), \\
&\!\!=\!\!&-\rs+j \qquad (\max(0,a-b)\le j\le a-1). 
\en
\end{lem}
We remark that the exchange relation \eqref{fww2} holds true 
with $f_{\la,\mu}(u,v)$ replaced by $\tilde{f}_{\la,\mu}(u,v)$.

Returning to the general $\lambda$, 
denote by $\mu_1\ge\cdots\ge\mu_m$ its column lengths 
(hence the transposed diagram is $\la'=(\mu_1,\cdots,\mu_m)$). 
\begin{lem}\lb{lem:app24}
\be
W_{\la}(u)&=&\Bigl(\prod_{i=2}^{m}
\ma{u_{i-1}-u_{i}-\rs}^{-\mu_{i}+1}
\\
&&\times
\prod_{1\le i<j\le l}\!\!f_{(1^{\mu_i}),(1^{\mu_j})}(u_i,u_j)\cdot
W_{(1^{\mu_1})}(u_1)\cdots W_{(1^{\mu_m})}(u_m)
\Bigr)\Bigl|_{u_i=u-(i-1)\rs\atop 1\le i\le m}.
\en
\end{lem}
\proof
Let $\bar{\lambda}$ be the diagram obtained by removing the last column of 
$\lambda$ so that $\lambda'=(\bar{\lambda}',\mu_m)$. We show
\bea
W_\la(u)=
\ma{u-v-(m-1)\rs}^{-\mu_{m}+1}
f_{\bar{\lambda},(1^{\mu_m})}(u,v)
W_{\bar{\lambda}}(u) W_{(1^{\mu_m})}(v)
\Bigl|_{v=u-(m-1)\rs}
\lb{lem24}
\ena
by induction on $\mu_m$. 
The lemma follows by repeated use of this equation. 

If $\mu_m=1$, then \eqref{lem24} is immediate from the definition. 
Assuming the statement is true for $\mu_m$
($\mu_{m-1}\ge\mu_m+ 1\ge 2$) we consider \eqref{lem24} with 
$\lambda'=(\bar{\lambda}',\mu_m+1)$.  
We have
\bea
&&\ma{u-v-(m-1)\rs}^{-\mu_{m}}
f_{\bar{\lambda},(1^{\mu_m+1})}(u,v)
W_{\bar{\lambda}}(u) W_{(1^{\mu_m+1})}(v)
\nn\\
&&
=\ma{u-v-(m-1)\rs}^{-\mu_{m}+1}\ma{u-v'-\mu_m-(m-1)\rs}^{-1}
f_{\bar{\lambda},(1^{\mu_m})}(u,v)f_{\bar{\lambda},(1)}(u,v')
\nn\\
&&\quad\times
\ma{v-v'-\mu_m}f_{(1^\mu_m),(1)}(v,v')
W_{\bar{\lambda}}(u) W_{(1^{\mu_m})}(v)W_{(1)}(v')
\Bigl|_{v'=v-\mu_m}.
\lb{lem25}
\ena
Let us verify that the right hand side of \eqref{lem25} 
(before specialization $v'=v-\mu_m$) 
is regular at $v=u-(m-1)\rs, v'=u-\mu_m-(m-1)\rs$. 
{}From the induction hypothesis, 
\be
\ma{u-v-(m-1)\rs}^{-\mu_{m}+1}f_{\bar{\lambda},(1^{\mu_m})}(u,v)
W_{\bar{\lambda}}(u) W_{(1^{\mu_m})}(v)
\en
is regular at $v=u-(m-1)\rs$, and  
\be
\ma{v-v'-\mu_m}f_{(1^\mu_m),(1)}(v,v')
W_{(1^{\mu_m})}(v)W_{(1)}(v')
\en
is regular at $v'=v-\mu_m$. 
Finally Lemma \ref{lem:app23} implies that 
\be
\tilde{f}_{\bar{\lambda},(1)}(u,v')W_{\bar{\lambda}}(u) W_{(1)}(v') 
\en
is regular at $v'=u-\mu_m-(m-1)\rs$, and 
\be
&&\ma{u-v'-\mu_m-(m-1)\rs}^{-1}
\frac{f_{\bar{\lambda},(1)}(u,v')}
{\tilde{f}_{\bar{\lambda},(1)}(u,v')}
\\
&&=\ma{u-v'-\mu_m-(m-1)\rs}^{-1}
\prod_{1\le k\le m-1\atop 1\le j\le \mu_k}
\frac{\ma{u-v'-(j-1)-k\rs}}{\ma{u-v'-j-(k-1)\rs}}
\en
is also regular (since $\mu_{m-1}\ge \mu_m+1$).

We let $v=u-(m-1)\rs$ in \eqref{lem25} and change the order of  
specialization.  
Using the induction hypothesis 
for the diagram $\tilde{\lambda}'=(\bar{\lambda}',\mu_m)$,   
we obtain 
\be
&&\ma{u-v'-\mu_m-(m-1)\rs}^{-1}\ma{u-(m-1)\rs-v'-\mu_m}
\\
&&\quad\times 
f_{\bar{\lambda},(1)}(u,v')f_{(1^\mu_m),(1)}(u-(m-1)\rs,v')
W_{\tilde{\la}}(u)W_{(1)}(v')\Bigl|_{v'=u-\mu_m-(m-1)\rs}
\\
&&=f_{\tilde{\la},(1)}(u,v')W_{\tilde{\la}}(u)W_{(1)}(v')
\Bigl|_{v'=u-\mu_m-(m-1)\rs}
\\
&&=W_{\la}(u).
\en
\qed

\subsection{Rectangular diagrams}

For a rectangular Young diagram $\la=(m^a)$, 
we write $W_m^{(a)}(u)=W_{(m^a)}(u)$. 
The following relations may be viewed as an analog of the $T$-system 
for the transfer matrices discussed in \cite{KNS94}.

\begin{prop}
\bea
&&f_{(m^a),(m^a)}(u,v)W^{(a)}_m(u)W^{(a)}_m(v)\Bigl|_{v=u-\rs}
\nn\\
&&\qquad =(-1)^{a-1}
f_{((m+1)^a),((m-1)^a)}(u,v)W^{(a)}_{m+1}(u)
W^{(a)}_{m-1}(v)\Bigl|_{v=u-\rs},
\lb{T1}\\
&&
\tilde{f}_{(m^a),(m^a)}(u,v)W^{(a)}_m(u)W^{(a)}_m(v)\Bigl|_{v=u-1}
\nn\\
&&\qquad =(-1)^{m-1}C^{(a)}_m
\tilde{f}_{(m^{a+1}),(m^{a-1})}(u,v)W^{(a+1)}_{m}(u)
W^{(a-1)}_{m}(v)\Bigl|_{v=u-1},
\lb{T2}
\ena
where
\be
C^{(a)}_m=\prod_{1\le s,t\le m}
\frac{\ma{a-1-(s-t)\rs}}{\ma{a-(1+s-t)\rs}}.
\en
Both sides of \eqref{T1},\eqref{T2} are well defined. 
\end{prop}

We sketch below the proof of \eqref{T1}.
First we check the regularity of both sides at $v=u-\rs$. 
For the right hand side, this can be shown from 
Lemma \ref{lem:app221}.
For the left hand side, we use Lemma \ref{lem:app23} to find that 
\be
\tilde{f}_{(m^a),(m^a)}(u,v)W^{(a)}_m(u)W^{(a)}_m(v)
\en
has poles of order at most $2(m-1)$ at $u=v-\rs$. 
Since 
\be
\frac{f_{(m^a),(m^a)}(u,v)}{\tilde{f}_{(m^a),(m^a)}(u,v)}
=O([u-v+\rs]_x^{2(m-1)})
\qquad (u\rightarrow v-\rs),
\en
the desired regularity follows. 
In the same way (using Lemma \ref{lem:app221}) we see that 
\be
f_{(m^a),((m-1)^a)}(u,v)W^{(a)}_m(u)W^{(a)}_{m-1}(v)
\en
has poles of order at most $(a-1)$ at $u=v-\rs$. 

Consider the expression
\be
&&A\equiv f_{(m^a),((m-1)^a)}(u,u')f_{(m^a),(1^a)}(u,v)
f_{((m-1)^a),(1^a)}(u',v)
W^{(a)}_m(u)W^{(a)}_{m-1}(u')W^{(a)}_1(v).
\en
{}From the definition of $W_\lambda(u)$, 
we have 
\be
A&\!\!=\!\!&\ma{u'-v-(m-1)\rs}^{a-1}f_{(m^a),(m^a)}(u,u')
W^{(a)}_m(u)W^{(a)}_m(u')
\\
&&\quad+O(\ma{u'-v-(m-1)\rs}^{a})\qquad 
(v\rightarrow u'-(m-1)\rs),
\\
&\!\!=\!\!&\ma{u-v-m\rs}^{a-1}f_{((m+1)^a),((m-1)^a)}(u,u')
W^{(a)}_{m+1}(u)W^{(a)}_{m-1}(u')
\\
&&\quad
+O(\ma{u-v-m\rs}^{a})\qquad 
(v\rightarrow u-m\rs).
\en
Writing $y=u-v-m\rs$, $y'=u'-v-(m-1)\rs$ and multiplying 
both sides by $\ma{u-u'-\rs}^{a-1}$ we have the equality of the form
\be
\varphi(y,y')
&=&\ma{y}^{a-1}\ma{y-y'}^{a-1}\psi(y-y')+O(y^a)
\qquad (y\rightarrow 0),
\\
&=&\ma{y'}^{a-1}\ma{y-y'}^{a-1}\psi'(y-y')+O(y^{\prime a})
\qquad (y'\rightarrow 0),
\en
where $\varphi(y,y')$, $\psi(y-y')$ and $\psi'(y-y')$ are regular near 
$y=y'=0$. 
This implies that $(-1)^{a-1}\psi(0)=\psi'(0)$, 
and \eqref{T1} follows. 
 

\end{document}